\documentclass[reqno,final]{amsart}

\usepackage{caption}
\usepackage[ocgcolorlinks,allcolors={blue},breaklinks]{hyperref}
\usepackage{amsfonts,latexsym,graphicx,amsmath,amssymb,bbm,enumitem,url}
\usepackage{algorithm}
\usepackage{algpseudocode}
\makeindex

\usepackage[color,notref,notcite]{showkeys}

\definecolor{labelkey}{rgb}{0.6,0,1}
\makeatletter
\def\BState{\State\hskip-\ALG@thistlm}
\makeatother
\allowdisplaybreaks

\definecolor{labelkey}{rgb}{0.6,0,1}

\usepackage[normalem]{ulem}
\newcounter{corr}
\definecolor{violet}{rgb}{0.580,0.,0.827}
\newcommand{\corr}[3]{\typeout{Warning : a correction remains in page
\thepage}
				\stepcounter{corr}        
				{\color{blue}\ifmmode\text{\,\sout{\ensuremath{#1}}\,}\else\sout{#1}\fi}
        {\color{red}#2}
        {\color{violet} \fbox{\thecorr}#3}}

\newcounter{cst}
\catcode`\@=11
\def\ctel#1{C_{\refstepcounter{cst}\@bsphack
\protected@write\@auxout{}%
           {\string\newlabel{#1}{{\thecst}{\thepage}}}\thecst}}
\catcode`\@=12

\newcounter{cexp}
\catcode`\@=11
\def\terml#1{T_{\refstepcounter{cexp}\@bsphack
\protected@write\@auxout{}%
           {\string\newlabel{#1}{{\thecexp}{\thepage}}}\thecexp}}
\catcode`\@=12

\newcommand{\mathbi}[1]{{\boldsymbol #1}}

\newcommand{\eop}{{\unskip\nobreak\hfil\penalty50
           \hskip2em\hbox{}\nobreak\hfil\mbox{\rule{1ex}{1ex} \qquad}
   \parfillskip=0pt
   \finalhyphendemerits=0\par\medskip}}

\newtheorem{theorem}{Theorem}[section]
\newtheorem{remark}[theorem]{Remark}

\newtheorem{definition}[theorem]{Definition}

\definecolor{shadecolor}{gray}{0.92}
\definecolor{TFFrameColor}{gray}{0.92}
\definecolor{TFTitleColor}{rgb}{0,0,0}

\newcommand{\ba}{\begin{array}{llll}   }
\newcommand{\bac}{\begin{array}{c}}
\newcommand{\bari}{\begin{array}{r}}
\newcommand{\ea}{\end{array}}

\newcommand{\ban}{\begin{array}{llll}}
\newcommand{\ean}{\end{array}}

\newcommand{\be}{\begin{equation}}
\newcommand{\ee}{\end{equation}}

\newcommand{\beqsys }{\beqtab \left \{ \begin{array}{l}}
\newcommand{\eeqsys }{\end{array} \right . \eeqtab }

\newcommand{\benum}{\begin{enumerate}}
\newcommand{\eenum}{\end{enumerate}}

\newcommand{\beqtab}{\begin{eqnarray}} 
\newcommand{\eeqtab}{\end{eqnarray}}

\newcommand{\bfn}{\mathbf{n}}

\newcommand{\disc}{{\mathcal D}}
\newcommand{\discC}{{\mathcal C}}

\newcommand{\dt}{{\delta\!t}}
\newcommand{\dtDisc}{\dt^{(n+\frac{1}{2})}}

\renewcommand{\div}{{\mathop{\rm div}}}
\newcommand{\divg}{\div}

\newcommand{\edge}{\sigma}

\newcommand{\edges}{{\mathcal E}}              

\newcommand{\ICinterp}{\mathcal{I}}

\newcommand{\mesh}{{\mathcal M}}

\newcommand{\norm}[2]{\| #1 \|_{#2}}

\renewcommand{\O}{\Omega}

\newcommand{\R}{\mathbb R}

\newcommand{\x}{\mathbi{x}}

\newcommand{\y}{\mathbi{y}}

\renewcommand{\norm}[2]{\left\Vert#1\right\Vert_{#2}}

\def\PiDc{\Pi_{\discC}}

\usepackage{xparse}
\DeclareDocumentCommand{\RPiD}{ O{\disc} O{,0} }{\Pi_{#1}(X_{#1#2})}

\def\Fdof#1{{\bm{\mathcal{F}}(#1,\R)}}
\makeatletter
\def\Fdof{\@ifnextchar[{\@with}{\@without}}
\def\@with[#1]#2{{\bm{\mathcal{F}}(#2;#1)}}
\def\@without#1{{\bm{\mathcal{F}}(#1,\R)}}
\makeatother

\newcommand{\Daru}{\mathbf{u}}

\def\RT0{\mathbb{RT}_0}

\begin{document}

	\title[]{An efficient implementation of mass conserving characteristic-based schemes in 2D and 3D}
	
		\author{Hanz Martin Cheng}
		\address{Department of Mathematics and Computer Science, Eindhoven University of Technology, P.O. Box 513, 5600 MB Eindhoven, The Netherlands. \texttt{h.m.cheng@tue.nl} \and
			School of Mathematics, Monash University, Victoria 3800, Australia.}
		\author{J\'er\^ome Droniou}
		\address{School of Mathematics, Monash University, Clayton, Victoria 3800, Australia.\\
			\texttt{jerome.droniou@monash.edu}}
	
	\date{\today}
		\maketitle
	%
	%
	%

		\begin{abstract}
			In this paper, we develop the ball-approximated characteristics (B-char) method, which is an algorithm for efficiently implementing characteristic-based schemes in 2D and 3D. Core to the implementation of numerical schemes is the evaluation of integrals, which in the context of characteristic-based schemes with piecewise constant approximations boils down to computing the intersections between two regions. In the literature, these regions are approximated by polytopes (polygons in 2D and polyhedra in 3D) and, due to this, the implementation in 3D is nontrivial. The main novelty in this paper is the approximation of the regions by balls, whose intersections are much cheaper to compute than those of polytopes. Of course, balls cannot fully tessellate a region, and hence some mass may be lost. We perform some adjustments, and also solve an optimisation problem, in order to yield a scheme that is both locally and globally mass conserving. This algorithm can achieve results that are similar to those obtained from an implementation which uses polytopal intersections, with a much cheaper computational cost. 
		\end{abstract}

	\section{Models and assumptions} \label{sec:PureAdvectionModel}
	
	\subsection{Introduction}
	
	In this paper, we introduce an algorithm for implementing characteristic-based schemes for a pure advection model  
		\begin{equation}\label{eq:advection}
		\begin{cases}
		\phi \dfrac{\partial c}{\partial t} 
		+ \div(\Daru c) 
		= 0 &\qquad \mbox{ on } Q_T:=\O \times (0,T),\\
		c(\cdot,0)=c_{\rm ini}&\qquad\mbox{ on }\Omega.
		\end{cases}
	\end{equation}
		Here, $T>0$, $\O$ is a polytopal domain in $\R^d$ ($d\ge 1$), the porosity $\phi$, and the velocity $\Daru$ are given, with $\Daru\cdot\bfn=0$ on $\partial\O$. The unknown $c(\x,t)$ represents the amount of material (a fraction) present at $(\x,t)$.  Note that the boundary is non-characteristic due to the assumption $\Daru\cdot\bfn=0$ on $\partial\O$, and thus no boundary conditions need to be enforced in \eqref{eq:advection}.

	The need to solve advection equations of the form \eqref{eq:advection} usually forms part of an operator splitting technique used to solve an advection-diffusion equation
	\begin{equation} \label{eq:advection-diffusion}
	\left\{\begin{array}{ll}
	\phi \dfrac{\partial c}{\partial t} 
	+ \div(\Daru c -\Lambda\nabla c) 
	= f(c) &\qquad \mbox{ on } Q_T:=\O \times (0,T),\\
	\Lambda\nabla c\cdot\bfn=0&\qquad\mbox{ on }\partial\O\times (0,T),\\
	c(\cdot,0)=c_{\rm ini}&\qquad\mbox{ on }\Omega,
	\end{array}\right.
	\end{equation}
	 where the source term $f$ and the diffusion tensor $\Lambda$ are given. Although \eqref{eq:advection} is presented with
	 $f=0$, it is not difficult to handle equations with nonzero reaction/source terms. This can be done, e.g., by a  splitting approach: solve the pure advection equation first (using the method presented here), then solve the reaction equation.
 
	 Advection-diffusion equations of the form \eqref{eq:advection-diffusion} are usually encountered in mathematical models for porous media flow (e.g. reservoir simulation, nuclear waste storage) \cite{E83-Mathematics-Reservoir,PR62}, and computational fluid dynamics (e.g. Navier-Stokes equations) \cite{NSbook89}, and are usually advection dominated. The diffusive component of the model is discretised separately, by mixed finite elements (MFEM), finite volumes, or other schemes that fit in the framework of the gradient discretisation method (GDM) \cite{gdm}, and will not be detailed in this paper. Here, we only focus on the implementation of characteristic-based schemes in \eqref{eq:advection}, such as the Eulerian Lagrangian Localised Adjoint Method (ELLAM) and the Modified Method of Characteristics (MMOC). The advantage of these schemes stems from the fact that they are based on characteristic methods, and thus capture the advective component of the PDE better than standard upwind schemes.	
	Several variants of the ELLAM, some of which are the finite element (FE) ELLAM \cite{CRHE90-ELLAM-main} and the finite volume (FV) ELLAM \cite{HR98-FV-ELLAM}, as well as a summary of their properties, have been presented in \cite{RC-02-overview}. One of the major issues faced when implementing characteristic-based schemes is the conservation of mass (both local and global). In order to achieve global mass conservation, some adjustments were performed on the MMOC, leading to the development of MMOC with adjusted advection (MMOCAA) \cite{DFP97-MMOCAA-main}. Although the MMOCAA achieves global mass conservation, it does not achieve local mass conservation. On the other hand, from its formulation, ELLAM satisfies global mass conservation; more recent variants of the ELLAM, such as the volume corrected characteristics mixed method (VCCMM) \cite{AH06,AH10-Fully-Conservative-ELLAM,AW11-stability-monotonicity-implementation}, achieve local volume conservation by adjusting the points tracked through the characteristics. These points may also be adjusted by following the algorithm proposed in \cite{D16-Opti-meshCorr}. Another way to achieve local volume conservation for characteristic-based schemes has been proposed in \cite{CDL18-GEM}. This is particularly useful for schemes with piecewise constant approximations, such as hybrid and mixed finite volume type schemes \cite{dro-10-uni}. As an example, in \cite{CDL18-GEM}, it was used to perform adjustments to make the HMM--ELLAM schemes in \cite{CD17-HMM-ELLAM-complete} locally mass conserving. More details about the convergence analysis and implementation of GDM characteristic-based schemes for \eqref{eq:advection-diffusion} and its applications to flows in porous media, can be found in \cite{CDL18-GEM,CDL17-convergence-ELLAM}. For schemes with piecewise constant approximations, evaluating the integrals arising from the discretisation of \eqref{eq:advection} boils down to computing intersections between polytopal regions (polygons in 2D and polyhedra in 3D). Although several algorithms are available for taking the intersection of polygons in 2D, they are quite expensive to implement in practice. Moreover, even though these methods can theoretically be extended to 3D, the main difficulty for a 3D implementation would come from taking intersections between polyhedra. Most of the polyhedral intersection algorithms in 3D are able to compute the intersection between two convex polyhedra efficiently, as in \cite{C16-convex-int,C92-convex-int,HMMN84-spacesweep-int,M78-int-convex}. However, even though the cells are initially convex, the tracked cell may not be convex. To our knowledge, the intersection of a convex polyhedron with a general polyhedron has only been dealt with in \cite{DMY93-convex-gen-int}, and even here, the computation of the intersection is not trivial or easy to implement.
	
	The purpose of this paper is to develop a feasible method to implement charac\-teristic-based schemes in 2D and 3D, whilst preserving the important properties of local and global mass conservation. The novelty of this paper is the idea of approximating the polytopal regions by balls (circles in 2D, spheres in 3D). By doing so, we convert the problem of computing polytopal intersections into that of computing intersections of balls, which is trivial to implement and has an essentially zero computational cost. Naively doing so will lead to a loss of mass, and hence we propose an adjustment algorithm which will help reduce the errors induced by this loss of mass. We then design to solve an optimisation problem, with both global and local mass conservations as constraints. We call this process the Ball-approximated characteristics (B-char) method.  Due to its formulation, the B-char method will yield a scheme that is both locally and globally mass conserving. 
	
	The paper is organised as follows. We start by giving some details on the assumptions on the data for the advection equation \eqref{eq:advection}. After which, we give a short summary of the ELLAM scheme used to discretise this equation in Section \ref{sec: ELLAM}. We also enumerate some of its mass conservation properties, and give a physical interpretation of the scheme. We then give a brief summary of how the ELLAM type schemes were implemented in the literature. The B-char method is introduced in Section \ref{sec:B-char}. In Section \ref{sec:numTests}, numerical tests are first performed in 2D in order to compare the performance of the B-char method with the ELLAM scheme obtained from polygonal intersections, with volume adjustments as described in \cite{CDL18-GEM}. In these tests, we see that the B-char method yields very similar results to the polygonal intersections, with a much cheaper computational cost. The B-char method is applied on some benchmark test cases; here too the results demonstrate the accuracy of the method. Finally, numerical tests are performed to show the applicability of the B-char ELLAM in 3D.


	\subsection{Assumptions on the data, and numerical setting} \label{sec:data-and-num-setting}
	We start by forming a mesh, i.e. a partition of $\O$ into polygonal (in 2D) or polyhedral (in 3D) sets. Following the notations in \cite[Definition 7.2]{gdm}, we then denote $\mathcal T=(\mesh, \edges)$ to be the set of cells $K$ and faces (edges in 2D) $\edge$ of our mesh, respectively. We also use $|K|$ to denote the volume (area in 2D) of a cell $K$. Throughout the article we assume the following properties:
	\begin{equation}\label{assump.global}
		\begin{aligned}
		&c_{\rm ini}\in L^\infty(\O),\\
		&\phi\in L^\infty(\O) \mbox{ is piecewise constant on } \mesh, \mbox{ and }\\
		&\mbox{there exists $\phi_*>0$ s.t. } \phi\ge \phi_*\mbox{ a.e.~on $\O$}.
		\end{aligned}
	\end{equation}
	
	Assumption \eqref{assump.global} simply states that the initial concentration inside the medium is bounded and  that the porosity $\phi$ of the medium does not vanish, which is natural in physical applications. The piecewise constant assumption on $\phi$ is also satisfied in practical applications; the value of $\phi$ on a cell $K$ will be denoted by $\phi_K$. As in \cite{CDL18-GEM}, we describe the numerical method in a general setting, to ensure that our algorithm applies at once to various possible spatial discretisations for the diffusion terms in \eqref{eq:advection-diffusion}. These can be dealt with using the GDM as shown in \cite{CDL18-GEM,CDL17-convergence-ELLAM}, and will not be discussed in further detail for this paper.
	We replace, in the weak formulation of the model, the continuous (infinite-dimensional) spaces and corresponding operators by a discrete (finite-dimensional) space and function reconstructions.  We then define a space-time discretisation $\discC=(X_\discC,\Pi_\discC,\ICinterp_{\discC},(t^{(n)})_{n=0,\dots, N})$, where
	\begin{itemize}
		\item $X_{\discC}$ is a finite-dimensional real space, describing the unknowns of the
		chosen scheme,
		\item $\Pi_\discC:X_\discC\to \mathbb{P}^0(\mesh)$ is a linear operator that reconstructs a piecewise constant function on the mesh $\mesh$
		from the unknowns,
		\item $\ICinterp_\discC$ is a rule to map $c_{\rm ini}$ onto an element
		$\ICinterp_\discC c_{\rm ini} \in X_{\discC}$,
		\item $0=t^{(0)}<t^{(1)}<\dots <t^{(N)} = T$ are the time steps, and we let
		$\dt^{(n+\frac12)}=t^{(n+1)}-t^{(n)}$.
	\end{itemize}
	Different choices of $\discC$ lead to different schemes (e.g. finite volume based methods, including hybrid ones with face unknowns like HMM \cite{dro-10-uni}, or mass-lumped finite element methods \cite{T06-FEbook}).
	
	Finally, we assume that 
	\be\label{hyp:velocity}
	\begin{aligned}
		&\Daru \in L^\infty(0,T;L^2(\O)^d) \mbox{ and }\div \Daru\in L^\infty(Q_T),
	\end{aligned}
	\ee
	and that $\Daru$ is approximated on each time interval $(t^{(n)},t^{(n+1)})$ by a function 
	\begin{equation}\label{hyp:un}
	\Daru^{(n+1)}\in H(\divg,\O) \mbox{ which is piecewise polynomial on } \mesh.
	\end{equation}
	
		\begin{remark}[Approximation of the velocity field]
		Although $\Daru$ is given in \eqref{eq:advection}, we use an approximation for the velocity field $\Daru$ in order to include the more general case where $\Daru$ comes from solving a PDE coupled to \eqref{eq:advection-diffusion}. For example, for flows in porous media, $\Daru$ usually comes from Darcy's law: Given a source term $g, \Daru$ should satisfy the PDE $-\divg(\Daru) = g$ on $\O $, with suitable boundary conditions.
	\end{remark}

	In the rest of the paper, the variables are only made explicit in the integrals when there is a risk of confusion. Otherwise we simply write, e.g., $\int_\O \phi d\x$.
	
	\section{ELLAM scheme for the advection--reaction equation} \label{sec: ELLAM}
	
 We multiply \eqref{eq:advection} with a sufficiently smooth function $\psi$, and perform integrations by parts. Using the identity
	\[
	\psi  \dfrac{\partial c}{\partial t}  =  \dfrac{\partial (c\psi)}{\partial t} - c\dfrac{\partial \psi}{\partial t},
	\]
    \eqref{eq:advection} gives, for any time interval $(t^{(n)},t^{(n+1)})$,
	\begin{align*}
	\int_{t^{(n)}}^{t^{(n+1)}}&\int_\O\phi(\x)\dfrac{\partial (c\psi)}{\partial t} (\x,t) d\x dt\nonumber\\
	&-
	\int_{t^{(n)}}^{t^{(n+1)}}\int_\O c(\x,t)
	\bigg[
	\phi(\x)\dfrac{\partial \psi}{\partial t} (\x,t) + \Daru(\x,t)\cdot \nabla \psi(\x,t)
	\bigg] d\x dt\nonumber = 0.\\
	\end{align*}
	To simplify the second term on the left hand side of the above equation,
	the ELLAM requires that test functions $\psi$ satisfy
	\begin{equation}\label{eq:testfunc}
	\phi\dfrac{\partial \psi}{\partial t} + \Daru\cdot \nabla \psi= 0 
	\quad\text{ on }\O \times(t^{(n)},t^{(n+1)}),
	\end{equation}
	with $\psi(\cdot,t^{(n+1)})$ given.
	The advection equation~\eqref{eq:advection} then leads to the relation
	\begin{equation} \label{eq:advection-weak}
	\int_\O\phi(\x)(c\psi) (\x,t^{(n+1)}) d\x
	-
	\int_\O\phi(\x)(c\psi) (\x,t^{(n)}) d\x
	= 0.
	\end{equation}
	
	We now write the ELLAM scheme, which consists of writing \eqref{eq:advection-weak} in the discrete context, 
	in which trial and test functions are replaced by reconstructions $\Pi_\discC$ applied to trial and
	test vectors in $X_\discC$.
	\begin{definition}[ELLAM scheme]\label{def:ELLAM}
		Given a space-time discretisation $\discC$, the ELLAM scheme for \eqref{eq:advection} reads as: find 
		$(c^{(n)})_{n=0,\ldots,N}\in X_{\discC}^{N+1}$ such that $c^{(0)}=\ICinterp_{\discC} c_{\rm ini}$
		and, for all $n=0,\ldots,N-1$, 
		$c^{(n+1)}$ satisfies
		\begin{equation}\label{conc-ELLAM}
		\int_{\O} \phi \PiDc c^{(n+1)} \PiDc z -\int_{\O} \phi \PiDc c^{(n)} \psi_z(t^{(n)}) = 0 \qquad \forall z \in X_{\discC}, 
		\end{equation}
		where $\psi_z$ is the solution to 
		\begin{equation}\label{Advection-test}
		\phi \partial_t \psi_z+\Daru^{(n+1)} \cdot \nabla \psi_z =0 \quad \text{ on } (t^{(n)},t^{(n+1)})\,,\mbox{ with $\psi_z(\cdot,t^{(n+1)})=\PiDc z$ }.
		\end{equation} 	
	\end{definition}
	\noindent Define the flow $F_{t}:\O\to\O$ such that, for a.e.\ $\x\in\O$,
	\begin{equation} \label{charac}
	\dfrac{dF_{t}(\x)}{dt} =\dfrac{\Daru^{(n+1)}(F_{t}(\x))}{\phi(F_t(\x))}\quad\mbox{ for $t\in [-T,T]$}, \qquad F_{0}(\x)= \x.
	\end{equation}
	Under Assumptions \eqref{assump.global} and \eqref{hyp:un}, the existence of this flow is proved in \cite[Lemma 5.1]{CDL17-convergence-ELLAM}.
	The solution to \eqref{Advection-test} is then understood in the sense:
	for $t\in (t^{(n)},t^{(n+1)}]$ and a.e.\ $\x\in\O$, $\psi_z(\x,t)=\Pi_\discC z(F_{t^{(n+1)}-t}(\x))$. In particular,
	\begin{equation}\label{eq:testfunc.ELLAM}
	\psi_z(\cdot,t^{(n)})=\PiDc z(F_{\dt^{(n+\frac12)}}(\cdot)).
	\end{equation}
	
	The construction of the B-char method in Section \ref{sec:B-char} will draw inspiration from a physical interpretation of the ELLAM, which uses the fact that $\Pi_\discC$ is a piecewise-constant reconstruction on a given mesh $\mesh$. For each cell $K\in\mesh$, we assume that there is $z_K\in X_\discC$ such that $\Pi_\discC z_K=\mathbbm{1}_K$, where $\mathbbm{1}_K$ is the function that has a value of 1 in $K$, and 0 elsewhere. Writing $\Pi_{\discC} c^{(k)}= \sum_{M\in\mesh} c_M^{(k)} \mathbbm{1}_M$ and taking $z_K$ as test function, \eqref{conc-ELLAM} and \eqref{eq:testfunc.ELLAM} give
	\[
	\int_{K} \phi \PiDc c^{(n+1)} d\x={}\int_{\O} \phi \sum_{M\in\mesh} c_M^{(n)} \mathbbm{1}_M(\x) \mathbbm{1}_K(F_{\dtDisc}(\x))d\x,
	\]
	which reduces to
	\begin{equation}\label{eq:phys-inter-ELLAM}
	|K|_\phi c_K^{(n+1)}= \sum_{M\in\mesh}|M\cap F_{-\dtDisc}(K) |_{\phi}c_M^{(n)},
	\end{equation}
	where $|E|_{\phi}=\int_E \phi$ is the available porous volume in a set $E\subset \R^d$.
	The term on the right hand side of \eqref{eq:phys-inter-ELLAM} tells us that the amount of material $c_K^{(n+1)}$ present in a particular cell $K\in\mesh$ at time $t^{(n+1)}$ is obtained by intersecting a tracked cell $F_{-\dtDisc}(K)$ and a residing cell $M$. This intersection can be interpreted as locating where the material in cell $K$ comes from, hence back-tracking the cell $K$ to $F_{-\dtDisc}(K)$, measuring which fraction of the material $c^{(n)}_M$ is taken from each $M\in\mesh$ (by taking their intersection), and depositing this fraction into the cell $K$.
	
	\subsection{Global mass conservation}
	Since the advection equation \eqref{eq:advection} usually comes from solving a model in computational fluid dynamics or engineering, we would want our numerical scheme to conserve global mass. Essentially, we would want an equation which tells us that  the change in $c$ is dictated by the amount of inflow/outflow and by the source term. In this case, due to the no-flow boundary conditions and the absence of a source term, this simply means that the amount of substance present at time $t^{(n+1)}$ should be the same as the amount of substance present at time $t^{(n)}$. The desired equation is thus given by
	\begin{equation} \label{mass-bal}
	\int_\O\phi(\x) c(\x,t^{(n+1)})d\x =\int_\O \phi(\x)c(\x,t^{(n)})d\x.
	\end{equation}
	
	It can easily be checked that the ELLAM scheme satisfies this property. Indeed, taking the sum over all $K\in\mesh$ in \eqref{eq:phys-inter-ELLAM} yields
	\[
	\sum_{K\in\mesh} c_{K}^{(n+1)}|K|_\phi = \sum_{M\in\mesh} c_M^{(n)} |M|_\phi,
	\]
	which is the discrete form of \eqref{mass-bal}.
	\begin{remark}[Achieving global mass conservation]
		We note here that the ELLAM scheme achieves global mass conservation due to
		\begin{equation} \label{eq:globalmassconsreq}
		\sum_{K\in\mesh} |M \cap F_{-\dtDisc}(K) |_\phi = |M|_\phi,
		\end{equation}
		 for all $M\in\mesh$. An analogue of this identity will be needed to ensure that the B-char method in Section \ref{sec:B-char} also achieves global mass conservation.
		\end{remark}
	
	\subsection{Local mass conservation}
	
\begin{figure}[h]
		\centering
		\begin{tabular}{c@{\hspace*{2em}}c}
			\includegraphics[width=0.4\textwidth]{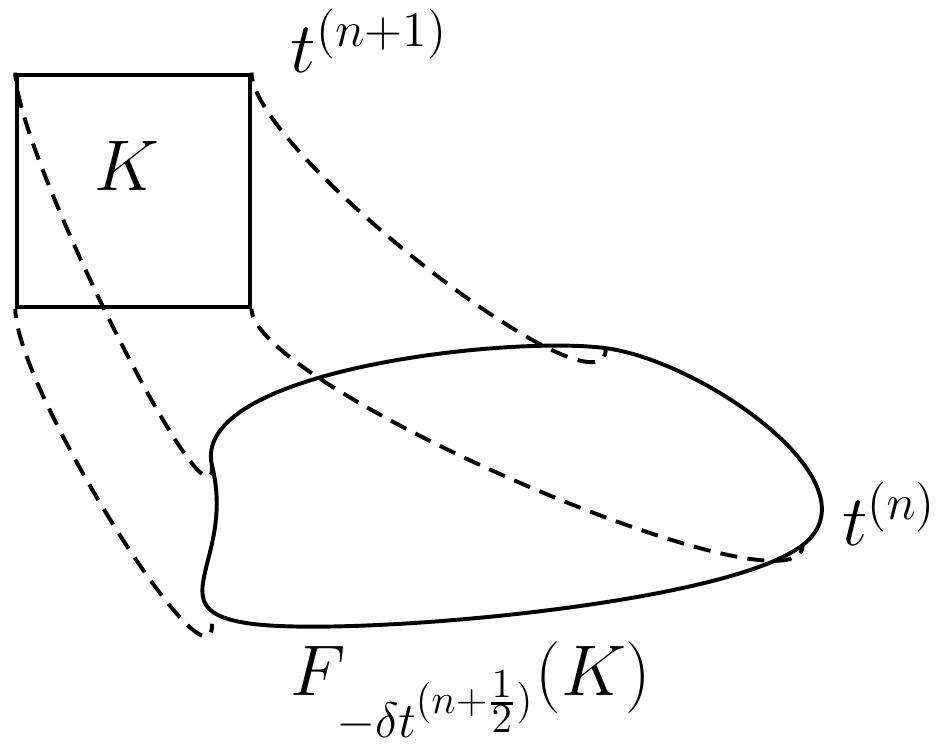} & 		\includegraphics[width=0.4\textwidth]{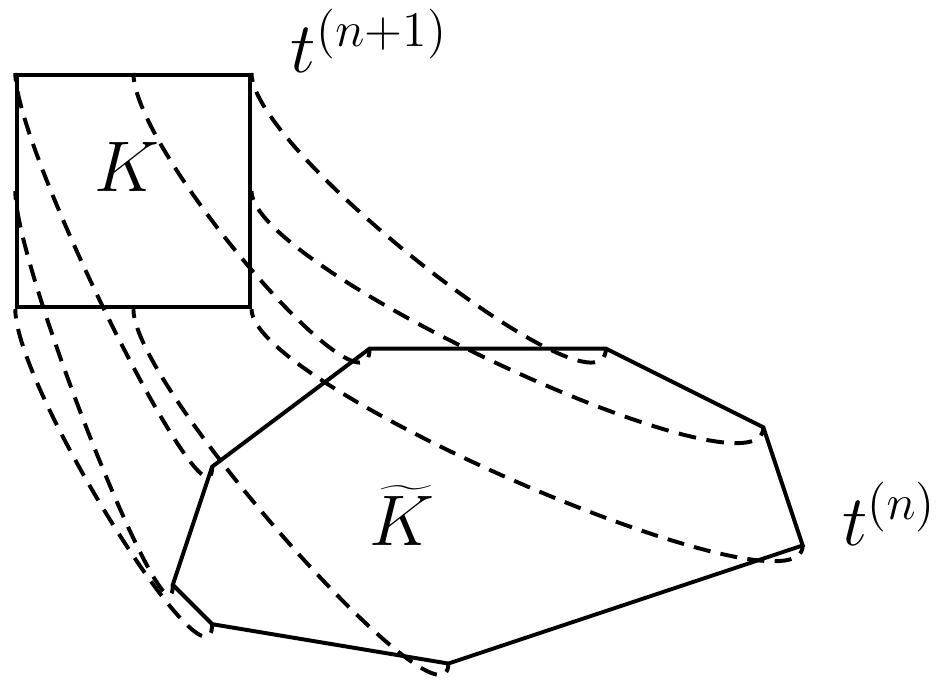}\\
		\end{tabular}
		\caption{ Trace-back region $\widehat{K}$ (left: exact; right: polygonal approximation $\widetilde{K}$.}
		\label{trace-back-regions}
	\end{figure}

	One of the main difficulties of implementing an ELLAM type scheme is the evaluation of the integral $\int_{\widehat{K}} \phi \PiDc c^{(n)} d\x$ for each cell $K$, where $\widehat{K}=F_{-\dtDisc}(K)$. In general, the region $\widehat{K}$ (see Figure \ref{trace-back-regions}, left) cannot be exactly described and hence, in the literature, it was approximated by polygons obtained from back-tracking the vertices, together with a number of points along the edges of the cell $K$. Figure \ref{trace-back-regions} (right) gives an illustration of the approximate trace-back region $\widetilde{K}$ obtained by tracking the vertices, together with the edge midpoints of the cell $K$.
	\begin{remark}[Reconstruction of polytopes]\label{rem:trackedPolyRecon}
		In 2D, most of the time, we can reconstruct the polygons approximating the trace-back region by following the tracked points in the same order as the original points, since it gives a well-defined polygon. However, in 3D, a face that is tracked may no longer be planar, and the original polyhedron faces need to actually be triangulated to ensure that a polyhedron is created after tracking.
		\end{remark}

	In general, $|\widetilde{K}|_\phi \neq |F_{-\dtDisc}(K)|_\phi$. However, the equality of these volumes is essential, otherwise the numerical scheme will not be able to preserve even a constant solution. Consider, for example, the simple case of a divergence free velocity field in \eqref{eq:advection}, with $\phi=1$ and $c_{\rm{ini}} = 1$. In this test case, the exact solution is given by $c(\x,t)=1$. In theory, upon implementing an ELLAM scheme with piecewise constant approximations for the unknown $c$, we should have the following simplified form of \eqref{eq:phys-inter-ELLAM} at the first time step:
	\begin{equation}\nonumber
	\begin{aligned}
	|K| c_K^{(1)} &= \sum_{M\in\mesh}|M\cap F_{-\frac{1}{2}}(K) |\times 1\\
	&= |F_{-\frac{1}{2}}(K)| \\ 
	&= |K|\qquad \mbox{(since $\mathbf{u}$ is divergence free)}.
	\end{aligned}
	\end{equation}
	However, due to the approximation of the trace-back region, we only have 
	\[
	|K| c_K^{(1)} = \sum_{M\in\mesh}|M\cap \widetilde{K} |\times 1 = |\widetilde{K}|
	\]
	and thus 
	\[
	c_K^{(1)} = \dfrac{|\widetilde{K}|}{|K|}\neq 1.
	\]
 This example shows that an inaccurate approximation of the volume of the tracked cell renders the numerical scheme unable to recover constant solutions. Hence, we need to perform some adjustments on the polygonal region $\widetilde{K}$ in order to yield $|\widetilde{K}|_\phi = |F_{-\dtDisc}(K)|_\phi$, which we shall define as the \textit{local volume constraint} for $K$. Several adjustment strategies which would lead to local mass conservation have been studied, as in \cite{AH06,CDL18-GEM,D16-Opti-meshCorr}. In particular, for local mass conservation to be achieved, we should have, for all $K\in\mesh$,
 \begin{equation} \label{eq:localmassconsreq}
  \sum_{M\in\mesh} |M\cap \widetilde{K}|_\phi = |F_{-\dtDisc}(K)|_\phi.
 \end{equation}
 
   For simplicity of exposition, we consider solenoidal fields, so that $\divg \mathbf{u} = 0$. The generalised Liouville's formula \cite[equation (26)]{CDL17-convergence-ELLAM} determines the evolution of available porous volume in a given domain: for any measurable set $A\subset \O$,
 	\begin{equation}\label{eq:genLiouville}
 	\dfrac{d}{dt}\int_{F_{t}(A)} \phi(\y)d\y = \int_{F_{t}(A)} \divg \mathbf{u}(\y)d\y.
 	\end{equation}	
Under the assumption that $\divg \mathbf{u} = 0$, \eqref{eq:genLiouville} with $A=K$  gives
 	\begin{equation}\label{eq:exactVol_traceback}
 	|F_{-\dtDisc}(K)|_\phi = |K|_\phi.
 	\end{equation}

	\section{B-char method}\label{sec:B-char}
	
	In this section, we present the idea of approximating the cells by balls, instead of the usual approximation using polygons. We will call this type of approximation the Ball-approximated Characteristics (B-Char). For each cell $K$, we choose $n_{K}$ points $C_{K,s}, (s=1,\dots,n_K)$ in its interior. We then assume that each of these points represents centers of disjoint balls $B_{K,s}$, with radius $r_{K,s}$, which are strictly inside cell $K$. The idea now is to distribute the porous volume in each cell $K$ over the balls $B_{K,s}$. To do so, we introduce a porous \emph{density} $\rho_K$ so that 
\begin{equation}\label{def:porous.density}
\rho_K \sum_{s=1}^{n_K}|B_{K,s}|_\phi=|K|_\phi.
\end{equation}
The porous density $\rho_K$ is fixed and does not change throughout the tracking.
 Here, the quantity $\rho_{K}|B_{K,s}|_\phi$ may be interpreted as an \emph{equivalent} porous volume inside the ball $B_{K,s}$.	The main interest of approximating the cells by balls is the fact that computing the intersection of balls is trivial compared to intersecting polytopes. As a consequence, the computational cost is greatly reduced. Moreover, this idea is easily applicable in both 2D and 3D.

	\begin{figure}[h]
		\centering
		\includegraphics[width=0.5\textwidth]{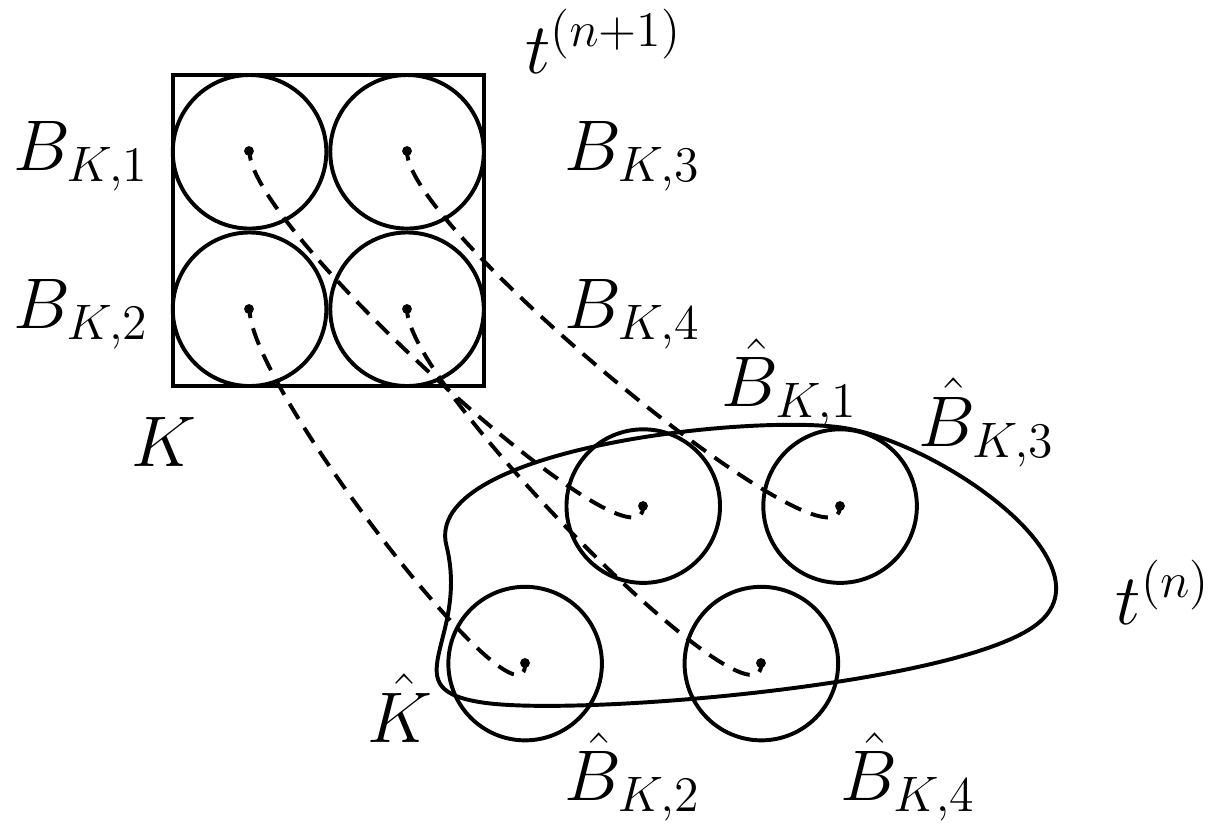}\\
		\caption{ Approximation of the trace-back region $\widehat{K}$ with balls.}
		\label{trace-back-ball}
	\end{figure}

	Upon working on the assumption that each ball, when tracked, remains a ball, the points $C_{K,s}$ are then tracked by solving \eqref{charac} to obtain $\widehat{C}_{K,s}$, which will be treated as the center of the tracked ball $\widehat{B}_{K,s}$ (see Figure \ref{trace-back-ball}). Of course, this assumption is not true in general, but gives a good approximation of the volumes, especially if the initial balls $B_{K,s}$ are not too large. Typically, a good enough approximation is obtained if the maximum radius of the balls $B_{K,s}$ is one-fourth or less of the diameter of the smallest face in the cell $K\in\mesh$.

			Since $\rho_K$ does not change throughout the tracking, the equivalent porous volume inside the tracked ball $\widehat{B}_{K,s}$ is given by $\rho_K|\widehat{B}_{K,s}|_\phi$. Denoting by $\widehat{\phi}_{K,s}$ the (unknown, at this stage) average porosity over the region covered by $\widehat{B}_{K,s}$, we may write
\begin{equation} \label{eq:porousBallVol}
\rho_K|\widehat{B}_{K,s}|_\phi = \rho_K \widehat{\phi}_{K,s} |\widehat{B}_{K,s}|.
\end{equation}
	Now, upon applying the generalised Liouville formula \eqref{eq:genLiouville} (recalling that $\divg\mathbf{u}=0$) with $A=B_{K,s}$, we have $|B_{K,s}|_\phi = |\widehat{B}_{K,s}|_\phi$,
which implies  
\begin{equation} \label{eq:expPhiB}
\phi_K |B_{K,s}| = \widehat{\phi}_{K,s} |\widehat{B}_{K,s}|.
\end{equation}

We now describe how to find the radius $\widehat{r}_{K,s}$ (and thus volume) of the tracked ball $\widehat{B}_{K,s}$, considering separately the case where $\phi$ is constant and where $\phi$ varies in the domain.

\begin{itemize}
\item \emph{$\phi$ constant in $\Omega$}. In this case, $\widehat{r}_{K,s}$ can be exactly computed. \eqref{eq:expPhiB} implies that the volume, and thus the radius, of the tracked ball is unchanged, i.e.
\begin{equation}\label{eq:ball_volChange}
\widehat{r}_{K,s} = r_{K,s}, \quad |\widehat{B}_{K,s}| = |B_{K,s}|. 
\end{equation} 
\item\emph{$\phi$ non-constant}. The radius $\widehat{r}_{K,s}$ cannot be exactly computed in general, only approximated with some additional computational cost. In this situation, aside from the center $C_{K,s}$, we also track points on the circumference of $B_{K,s}$. Considering the two-dimensional case for illustration, letting $C_{K,s}=(x_{K,s},y_{K,s})$, we also track the four cardinal points $C_{K,s_1}= (x_{K,s},y_{K,s}+r_{K,s})$, $C_{K,s_2}= (x_{K,s},y_{K,s}-r_{K,s})$, $C_{K,s_3}= (x_{K,s}+r_{K,s},y_{K,s})$, $C_{K,s_4}= (x_{K,s}-r_{K,s},y_{K,s})$.
The radius $\widehat{r}_{K,s}$ of the tracked ball $\widehat{B}_{K,s}$ is then approximated by
	\begin{equation}\label{eq:rTrackedBall}
	\widehat{r}_{K,s} \approx  \frac{1}{4}\sum_{j=1}^4\widehat{r}_{K,s_j},
	\end{equation} 
  where $(\widehat{r}_{K,s_j})_{j=1,\dots,4}$ denote the distances between the tracked center $\widehat{C}_{K,s}$ and the tracked cardinal point $(\widehat{C}_{K,s_j})_{j=1,\ldots,4}$. Of course, we may track more than 4 points on the circumference of the ball in order to get a better approximation of $\widehat{r}_{K,s}$. Finding the optimal number of points to be tracked in order to get a good approximation for $\widehat{r}_{K,s}$ depends on how strongly we expect the velocity field to distort the tracked region, and will be a topic for future research. The volume $|\widehat{B}_{K,s}|$ can then be obtained from $\widehat{r}_{K,s}$ defined by \eqref{eq:rTrackedBall}. 
\end{itemize}

 After computing $|\widehat{B}_{K,s}|$, we use \eqref{eq:expPhiB} to determine the value of $\widehat{\phi}_{K,s}$. Similarly to \eqref{def:porous.density}, we can then write the porous volume inside each tracked cell in the following form:
	\begin{equation}\label{eq:volTracked}
	|F_{-\dtDisc}(K)|_\phi = \sum_{s=1}^{n_K}  \rho_{K}\widehat{\phi}_{K,s}|\widehat{B}_{K,s}|.
	\end{equation}

	\subsection{Initial approximation for the volume of intersecting regions}
	We now describe the process for obtaining an initial approximation for the volume of the intersecting regions $|F_{-\dtDisc}(K) \cap M|_\phi \approx V_{\widehat{K},M}$.
	
	We start by recalling that $|F_{-\dtDisc}(K) \cap M|_\phi$ is interpreted as the amount of material in a cell $K$ that comes from a residing cell $M$. In the context of approximation by balls, this reads: each ball $B_{K,s}$ contains an amount of material from some residing balls $B_{M,m}$ transported by the flow, given by $\rho_M\phi_M |\widehat{B}_{K,s} \cap B_{M,m}|$. 
		
	\begin{remark}\label{rem:porousBallInt}
		The choice of $\rho_M\phi_M$ (instead of $\rho_K\phi_K$) in the initial approximation  $\rho_M\phi_M |\widehat{B}_{K,s} \cap B_{M,m}|$  comes from the interpretation that the amount of material present in $B_{K,s}$  is obtained by measuring how much of the material is taken from each $B_{M,m}$, and by depositing this material into the ball $B_{K,s}$. 
	\end{remark}
	
	An initial approach for approximating $V_{\widehat{K},M}$ would then involve taking the sum of the masses of the balls in a residing cell $M$, intersected with the tracked balls that originated from cell $K$, that is $\sum_{s=1}^{n_K}\sum_{m=1}^{n_M} \rho_M \phi_M |\widehat{B}_{K,s}\cap B_{M,m}|$. However, since there are gaps between the residing balls, $\sum_{M\in\mesh} \sum_{m=1}^{n_M} |\widehat{B}_{K,s}\cap B_{M,m}| \neq |\widehat{B}_{K,s}|$. This will lead to a loss in volume, which will in turn lead to a loss of mass conservation and a poor approximation. Instead, we use this to compute 	
	\begin{equation}\label{eq:fracMaterial}
		\dfrac{ \rho_M \phi_M|\widehat{B}_{K,s} \cap  B_{M,m}|}{\sum_{M\in\mesh} \sum_{\ell=1}^{n_{M}} \rho_M \phi_M|\widehat{B}_{K,s} \cap  B_{M,\ell}|},
	\end{equation}
	which represents the fraction of the mass in $\widehat{B}_{K,s}$ that comes from $B_{M,m}$. From this, we then see that
	\[
	\rho_K\widehat{\phi}_{K,s}|\widehat{B}_{K,s}|\dfrac{ \rho_M\phi_M|\widehat{B}_{K,s} \cap  B_{M,m}|}{\sum_{M\in\mesh} \sum_{\ell=1}^{n_{M}} \rho_M\phi_M|\widehat{B}_{K,s} \cap  B_{M,\ell}|}
	\]
	is the actual amount of mass in the tracked ball $\widehat{B}_{K,s}$ that comes from $B_{M,m}$. The quantity $V_{\widehat{K},M}$ is then computed by taking the sum over all tracked balls $\widehat{B}_{K,s}$ and residing balls $B_{M,m}$, given by 
	\begin{equation}\label{eq:ball_int}
	V_{\widehat{K},M} :=\sum_{s=1}^{n_{K}} \rho_K\widehat{\phi}_{K,s}|\widehat{B}_{K,s}| \dfrac{ \sum_{m=1}^{n_{M}}  \rho_M\phi_M |\widehat{B}_{K,s} \cap  B_{M,m}|}{\sum_{M\in\mesh} \sum_{\ell=1}^{n_{M}} \rho_M\phi_M |\widehat{B}_{K,s} \cap  B_{M,\ell}|}.
	\end{equation}

\subsection{Mass conservation for the B-char method}
	Since $V_{\widehat{K},M}$ are approximations to $|F_{-\dtDisc}(K) \cap M|_\phi$, in order to achieve local mass conservation, we should have an analogue of \eqref{eq:localmassconsreq}, given by $\sum_{M\in\mesh} V_{\widehat{K},M} = |F_{-\dtDisc}(K)|_\phi$. We can easily check that $V_{\widehat{K},M}$ in \eqref{eq:ball_int} satisfies this relation by using \eqref{eq:volTracked}. 
	Hence, the approximation of $|F_{-\dtDisc}(K) \cap M|$ by $V_{\widehat{K},M}$ in \eqref{eq:ball_int} leads to a scheme that is locally mass conserving. However, $\sum_{K\in\mesh} V_{\widehat{K},M} \neq |M|_\phi$,  which means that \eqref{eq:globalmassconsreq} is not satisfied, and thus global mass conservation is not achieved.   This indicates that some adjustments need to be performed on $V_{\widehat{K},M}$. 
	
		\begin{remark}[Overlapping balls] 
Two given tracked balls, say $\widehat{B}_{K,s_i}$ and $\widehat{B}_{L,s_j}$, may overlap. In this case, the mass in the overlapping region $\widehat{B}_{K,s_i} \cap \widehat{B}_{L,s_j}$ would be allocated twice (once each for the contribution to $B_{K,s_i}$ and $B_{L,s_j}$). Recalling Remark \ref{rem:porousBallInt}, this would imply that excess mass is deposited into cell $K$ or $L$ (or both), since these are the cells where the tracked balls $\widehat{B}_{K,s_i}$ and $\widehat{B}_{L,s_j}$ come from, respectively. This corresponds to a loss in mass that should be deposited in another cell, say $M$. This is implicitly reflected by the failure of the initial approximation $V_{\widehat{K},M}$ defined in \eqref{eq:ball_int} to achieve global mass conservation. The algorithm provided below to adjust $V_{\widehat{K},M}$, and recover both global and local mass conservations, also corrects the excessive and missing mass deposits induced by overlapping balls.
		\end{remark}
	
		Upon indexing each tracked cell $K_i$ and each residing cell $M_j$, $i,j = 1,\dots n_c$, let $A^{(n)}$ be the matrix with entries $a_{ij}^{(n)}$, where $a_{ij}^{(0)}= V_{\widehat{K}_i,M_j}$. In short, each entry $a_{ij}^{(n)}$ of the matrix $A^{(n)}$ gives an approximation of the volume $|F_{-\dtDisc}(K_i)\cap M_j|_\phi$. In terms of the matrix $A^{(n)}$, we would thus need 
	\begin{equation} \label{eq:localmassconsMat}
	\sum_{j=1}^{n_c} a_{ij}^{(n)} = |F_{-\dtDisc}(K_i)|_\phi \quad \mathrm{ for } \quad i=1,\dots,n_c
	\end{equation}
	in order to achieve \eqref{eq:localmassconsreq}, which will lead to local mass conservation. 
	Owing to \eqref{eq:exactVol_traceback}, the right hand side $|F_{-\dtDisc}(K_i)|_\phi$ of \eqref{eq:localmassconsMat} is equal to $|K_i|_\phi$, which can be easily and exactly computed since $\phi$ is constant in each cell. 
  Now, in order to have a globally mass conserving scheme, we would need to satisfy \eqref{eq:globalmassconsreq}, which in this context is equivalent to
	\begin{equation}\label{eq:globalmassconsMat}
	\sum_{i=1}^{n_c} a_{ij}^{(n)} = |M_j|_\phi \quad \mathrm{ for } \quad j=1,\dots,n_c.
	\end{equation}
	
\begin{figure}[h]
	\centering
	\includegraphics[width=0.5\textwidth]{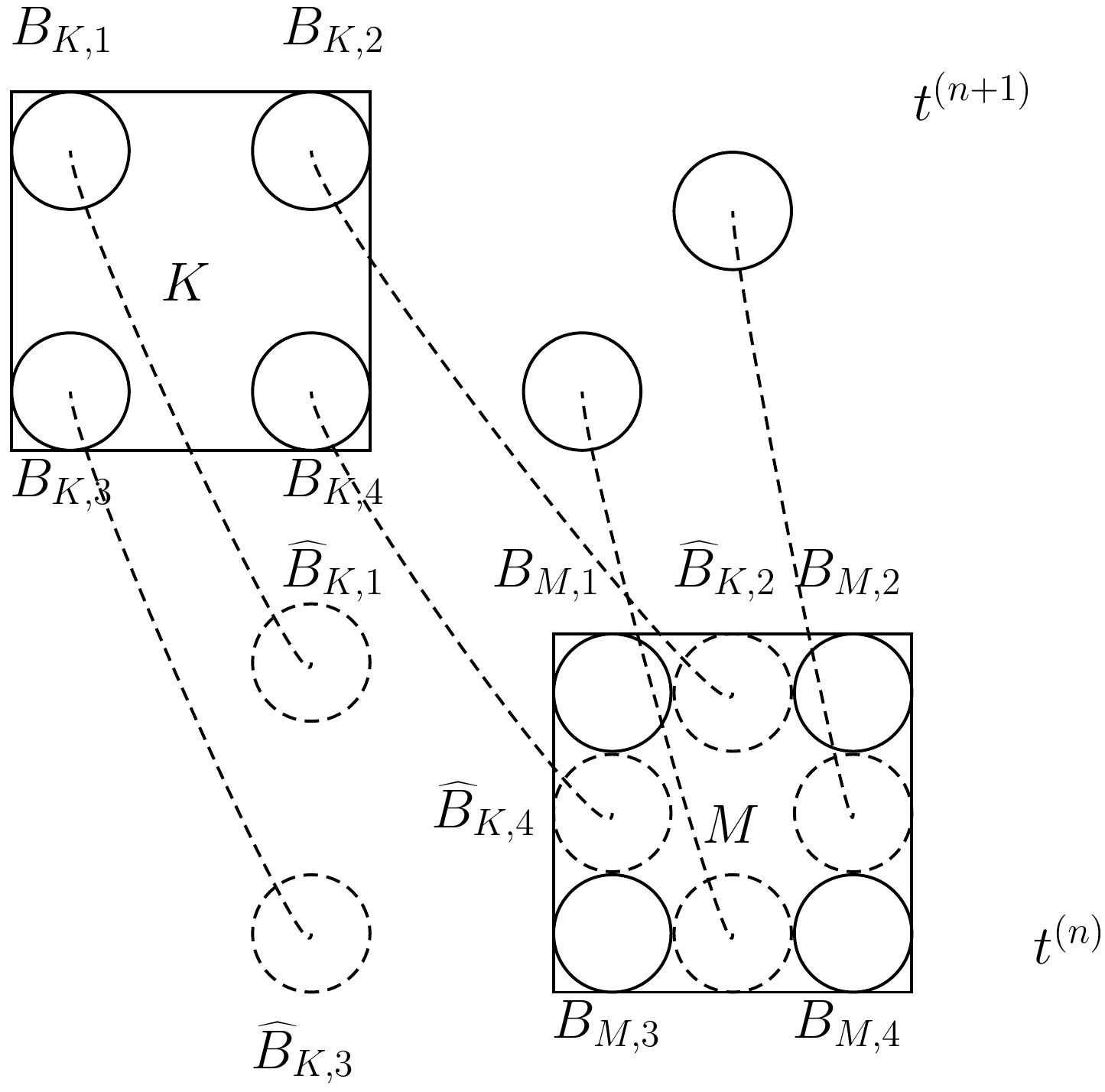}\\
	\caption{ Tracking into a void.}
	\label{trace-back-void}
\end{figure}

	To build a matrix which satisfies \eqref{eq:localmassconsMat}  and \eqref{eq:globalmassconsMat}, we start with the assumption that $\sum_{i=1}^{n_c} a_{ij}^{(0)} > 0 $ for $j=1,\dots,n_c$. This means that at least one of the balls $B_{M_j,\ell}$ of a residing cell $M_j$ intersects with at least one of the balls $B_{\widehat{K}_i,s}$ that has been tracked from $K_i$. The only time  $\sum_{i=1}^{n_c} a_{ij}^{(0)} = 0 $ holds is when the flow traces everything into voids. For example, in Figure \ref{trace-back-void}, none of the balls $\widehat{B}_{K,\ell}$ tracked from $K$, nor did any of the other balls that were tracked from time $t^{(n+1)}$, intersect with the residing balls $B_{M,\ell}$ in $M$, and so $\sum_{i=1}^{n_c} a_{ij}^{(0)} = 0 $ for the $j$ corresponding to the particular cell $M$. This can easily be resolved by increasing the number of balls, and making sure that each residing cell is tightly packed with balls. Another instance when $\sum_{i=1}^{n_c} a_{ij}^{(0)} = 0 $ is when there is a very strong inflow or outflow at the boundary of the domain, which is not the case for \eqref{eq:advection}.
\begin{remark}[Presence of inflow or outflow]
	 In the presence of inflows or outflows, we may avoid $\sum_{i=1}^{n_c} a_{ij}^{(0)} = 0 $ by one of the following options:
	\begin{enumerate}
		\item[a.] Take a smaller time step so that the region does not get emptied out.
		\item[b.] Create ghost cells at the boundary of the domain so that the cells that get emptied out are the ghost cells. 
	\end{enumerate} 
	\end{remark}
Now, start by setting 
	\[
	a_{ij}^{(n+\frac{1}{2})} = \frac{|M_j|_\phi}{\sum_{i=1}^{n_c} a_{ij}^{(n)}} a_{ij}^{(n)},
	\] 
	which yields
	\begin{equation} \nonumber
	\sum_{i=1}^{n_c} a_{ij}^{(n+\frac{1}{2})} = |M_j|_\phi.
	\end{equation}
	Next, we set 
	\begin{equation}\label{eq:volAdj}
	a_{ij}^{(n+1)} = \frac{|F_{-\dtDisc}(K_i)|_\phi}{\sum_{j=1}^{n_c} a_{ij}^{(n+\frac{1}{2})}} a_{ij}^{(n+\frac{1}{2})}
	= \frac{|F_{-\dtDisc}(K_i)|_\phi |M_j|_\phi}{\sum_{j=1}^{n_c} a_{ij}^{(n+\frac{1}{2})}\sum_{i=1}^{n_c} a_{ij}^{(n)}} a_{ij}^{(n)},
	\end{equation}
	so that
	\begin{equation} \nonumber
	\sum_{j=1}^{n_c} a_{ij}^{(n+1)} = |F_{-\dtDisc}(K_i)|_\phi.
	\end{equation}
	Essentially, the adjustments perform the following two-step process:
	\begin{itemize}
		\item  Firstly, we redistribute, according to a proportion evaluated by $a_{ij}^{(n)}$, the mass along each intersecting region so that global mass conservation \eqref{eq:globalmassconsMat} is achieved.
		\item Next, we redistribute the mass so that local mass conservation \eqref{eq:localmassconsMat} is achieved. 	
	\end{itemize}
	Intuitively, we can see that these adjustments involve re-distributing the errors and hence scaling them down in each iterate. A naive approach would involve iterating this process, stopping only when the error in global and local mass conservations are less than a certain tolerance value. However, there is no guarantee that such a result is achievable. A more efficient approach would involve, after taking $N$ iterations and arriving at the matrix $A^{(N)}$, solving a minimisation problem. In practice, we found that taking $N$ to be such that the maximum error in mass conservation is at most $5\%$ is sufficient to give a good initial approximation.
	
	We then assign an unknown corresponding to each entry of $A^{(N)}$, which gives us $n_c \times n_c$ unknowns. For $i,j$ such that $a_{ij}^{(N)} = 0$, the corresponding unknown is fixed to 0. This tells us that if no intersection has been detected between a tracked cell $\widehat{K}_i$ and a residing cell $M_j$, then our adjustment algorithm should not introduce any volume into these regions. Hence, the number of unknowns in our new system is equal to the number of nonzero entries $n_z$ in $A^{(N)}$. From the nonzero entries of $A^{(N)}$ construct a $2n_c \times n_z$ matrix $\hat{A}$ in the following manner: Write 
	\[
	A^{(N)} = \begin{bmatrix}
	\mathbf{r}_1  \\
	\mathbf{r}_2  \\
	\vdots \\
	\mathbf{r}_{n_c} 
	\end{bmatrix},
	\]
	and denote by $\hat{\mathbf{r}}_i$ ($i=1,\dots,n_c$) the row vector of size $\le n_c$ obtained by removing the zero entries in $\mathbf{r}_i$. The first $n_c$ rows of $\hat{A}$ are then formed by the $n_c\times n_z$ matrix
	\[  \begin{bmatrix}
	\hat{\mathbf{r}}_1 & \mathbf{0}  &  \cdots & \mathbf{0}\\
	\mathbf{0} & \hat{\mathbf{r}}_2 & \ddots & \mathbf{0}\\
	\vdots &  \ddots & \ddots  & \mathbf{0}\\
	\mathbf{0} & \mathbf{0} &  \mathbf{0} & \hat{\mathbf{r}}_{n_c}
	\end{bmatrix}.
	\]
	That is, we stagger the vectors $(\hat{\mathbf{r}}_i)_{i=1,\ldots,n_c}$ so that the coefficients of $\hat{\mathbf{r}}_{j+1}$ start at the column after the last coefficient of $\hat{\mathbf{r}}_j$, and we pad each row with zeros to ensure we obtain an $n_c \times n_z$ matrix.

	The latter $n_c$ rows of the matrix $\hat{A}$ are then formed in the following manner: for the $n_c+m$th row, we look for the nonzero entries $a_{j,m}$ ($j=1,\dots n_c$) in column $m$ of $A$. For each corresponding $j$, we then find the column $\hat{\jmath}$ corresponding to where the coefficient $a_{j,m}$ resides in the first $n_c$ rows of $\hat{A}$. We then set $\hat{a}_{n_c+m,\hat{\jmath}} = a_{j,m}$. As an example, if 
	\[
	A^{(N)} = \begin{bmatrix}
	a &  b & c & \cdots \\
	d & 0 & e & \cdots \\
	\vdots  &\vdots &\vdots &\ddots 
	\end{bmatrix},
	\] 
	then $\hat{A}$ would be (with the line separating the first $n_c$ rows from the last $n_c$ rows of $\hat{A}$)
	\[
	\hat{A}=	\begin{bmatrix}
	a & b & c & 0 & 0& 0 &\cdots \\
	0 & 0 & 0 & d & e & 0 &\cdots \\
	\vdots & \vdots & \vdots & \vdots & \vdots & \vdots& \vdots \\
	\hline
	a & 0 & 0 & d & 0 & \cdots& \cdots\\
	0 & b & 0 & 0 & 0 &\cdots& \cdots\\
	0 & 0 & c & 0 & e & \cdots& \cdots\\
		\vdots & \vdots & \vdots & \vdots & \vdots &\vdots & \ddots
	\end{bmatrix}.
	\] 
	In practice, the last $n_c$ rows of $\hat{A}$ are assembled simultaneously with its first $n_c$ rows. In the example above, after setting $a,b,c$ as the first, second and third entries of the first row, $a,b,c$ are simultaneously set to be the first, second and third entries of the $n_c+1, n_c+2,$ and $n_c+3$th row, respectively. After which, when $d$ and $e$ were set to be the fourth and fifth entry of the second row, they were also set to be the fourth and fifth entries of the $n_c+1$ and the $n_c+3$th row. In this manner, the whole matrix $\hat{A}$ is rather easy to assemble.
	
	The essential property of $\hat{A}$ is that when it is multiplied by an $n_z \times 1$ vector $\mathbf{1}$ consisting of all ones, we have
	\[
	\hat{A} \mathbf{1} = \begin{bmatrix}
	\sum_{j=1}^{n_c} a_{1j}^{(N)}\\
	\vdots\\
	\sum_{j=1}^{n_c} a_{n_cj}^{(N)}\\
	\sum_{i=1}^{n_c} a_{i1}^{(N)}\\
	\vdots \\
	\sum_{i=1}^{n_c} a_{in_c}^{(N)}\\
	\end{bmatrix}.
	\]
The sums in the right-hand side correspond to the quantities that must be fixed to certain values in order to achieve local and global mass balance, see \eqref{eq:localmassconsMat} and \eqref{eq:globalmassconsMat}. To obtain local and global mass conservation, we therefore solve the system 
	\begin{equation} \label{eq:opti_masscons}
	\widehat{A} (\mathbf{1}+\mathbf{x}) = \mathbf{b} := \begin{bmatrix} \mathbf{b}_\mathrm{loc}\\\mathbf{b}_\mathrm{glob} \end{bmatrix}, 
	\end{equation} 
	where $\mathbf{b}_\mathrm{loc}$ and $\mathbf{b}_\mathrm{glob}$ are the vectors containing the local and global mass constraints given by the right hand side of \eqref{eq:localmassconsMat} and \eqref{eq:globalmassconsMat}, respectively. Letting $\mathbf{x}=(x_j)_{j=1,\ldots,n_c}$ for all $j=1,\ldots,n_c$ we can view $1+x_j$ as an adjustment (in terms of scaling) of the approximations $a_{ij}^{(N)}$ of $|F_{-\dtDisc}(K_i)\cap M_j|_\phi$.

In general, a tracked cell intersects more than one residing cell, and so $n_z \geq 2n_c$, and hence the system is underdetermined and we have to select one of its solutions. This is done through the following minimisation problem: minimise $\mathbf{x}^T\mathbf{x}$ subject to the local and global mass constraints \eqref{eq:opti_masscons}. Moreover, since these quantities represent volumes, we also impose the constraint that each coefficient in $(\mathbf{1}+\mathbf{x})$ is positive. Finally, we impose that each entry of the vector $(\mathbf{1}+\mathbf{x})$ should be less than or equal to 2 (so that the maximum change is doubling a given volume). Essentially, this tells us that we want to achieve global and local mass conservation with minimal adjustments on the computed/approximated volumes, which makes sense since we assume that these intersections have been well-approximated. In terms of computational cost, solving the minimisation problem is not too expensive since a tracked cell usually only intersects a few residing cells (as long as the velocity field and the mesh are not too irregular), and hence the matrix $\hat{A}$ is usually sparse.

	\begin{remark}[Non-solenoidal fields]
		The B-char method may also be applied with non-solenoidal fields. In general, given a velocity field $\Daru$, the idea is to approximate the radius $\widehat{r}_{K,s}$, and hence the volume $|\widehat{B}_{K,s}|$, by \eqref{eq:rTrackedBall}. The quantity $\widehat{\phi}_{K,s}$ is then obtained from an analogue of \eqref{eq:expPhiB},   
		which comes from applying the generalised Liouville formula \eqref{eq:genLiouville}. 
		The quantity $\rho_K\widehat{\phi}_{K,s}$ and the tracked balls $\widehat{B}_{K,s}$ with radius $\widehat{r}_{K,s}$ can then be used in \eqref{eq:ball_int} for the initial approximation of $|F_{-\dtDisc}(K) \cap M|_\phi$. 
	\end{remark}

	\subsection{Summary of the B-char method}
	To summarise, the B-char method consists of approximating each cell $K\in\mesh$ by a collection of balls $B_{K,s}$ and their trace-back regions $F_{-\dtDisc}(K)$ by tracking back the centres of the balls $B_{K,s}$ to obtain the centres of the tracked balls $\widehat{B}_{K,s}$. We then perform the steps outlined in Algorithm \ref{alg:Vol_adj}:
	\begin{algorithm}
		\caption{Volume approximation and adjustment}\label{alg:Vol_adj}
		\begin{algorithmic}[1]
			\For {$i,j=1$ to $n_c$}
			
			\State Compute an initial approximation $V_{\widehat{K}_i,M_j} $ to $|F_{-\dtDisc}(K_i)\cap M_j|$ as in \eqref{eq:ball_int}.
			\EndFor
			\item Form a matrix $A^{(0)}$ with entries $a_{ij}^{(0)}=V_{\widehat{K}_i,M_j}$.
			\For{$n=1$ to $N$}
			\State Decrease the error in global mass conservation by updating the matrix $A^{(n)}$ with entries $a_{ij}^{(n)}$ as in \eqref{eq:volAdj}.
			\EndFor
			\item Find the minimal change in the approximated volumes $a_{ij}^{(N)}$ so that the constraints for global and local mass conservation \eqref{eq:opti_masscons} are satisfied. 
		\end{algorithmic}
	\end{algorithm}

	\section{Numerical tests} \label{sec:numTests}
\begin{figure}[h]
		\centering
		\includegraphics[width=0.8\textwidth]{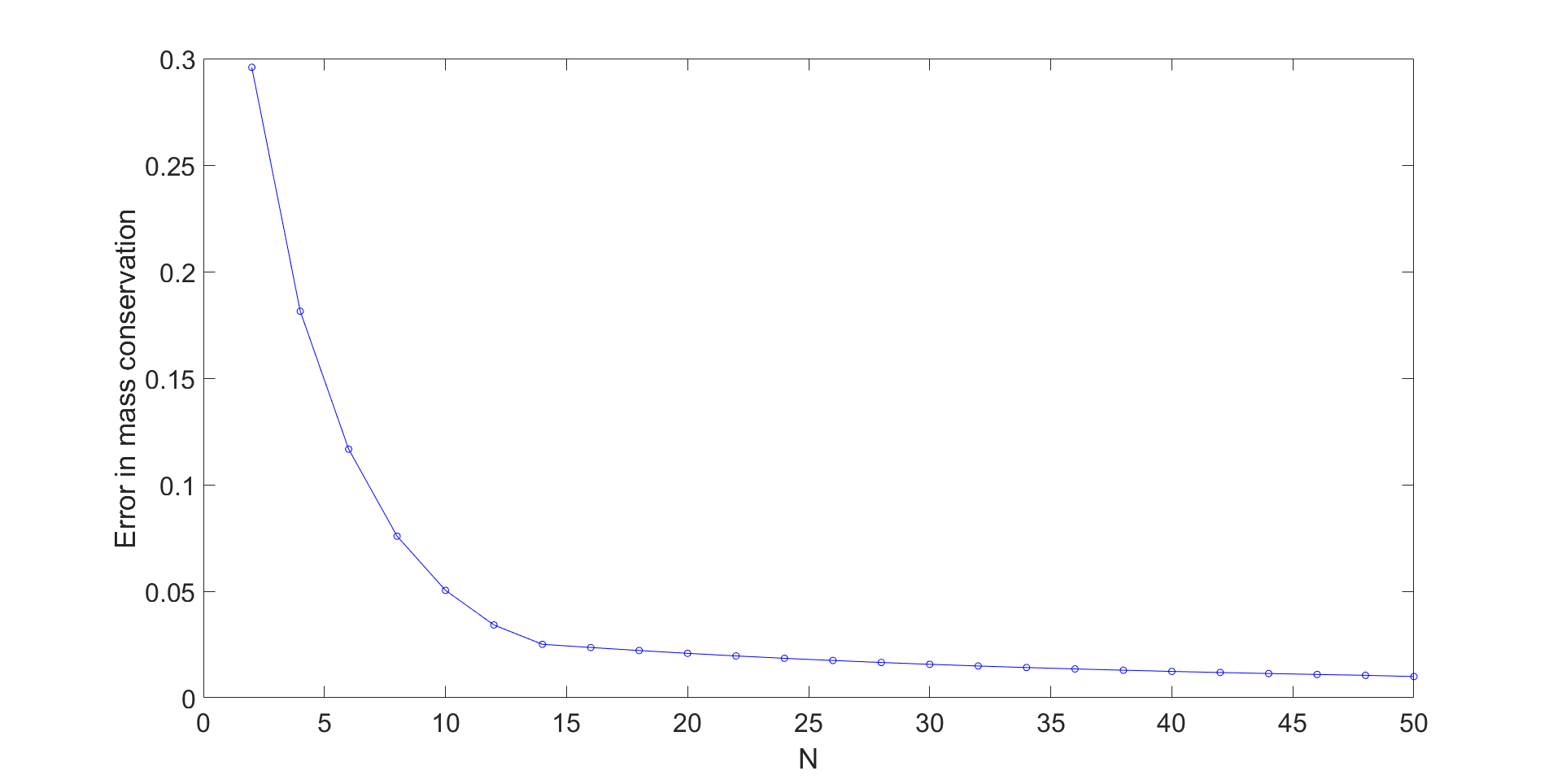} \\
		\caption{Maximum (relative) error in mass conservation for different choices of N.}
		\label{plot_N_massCons}
	\end{figure}

	In this section, we perform numerical tests on Cartesian type meshes for the pure advection equation \eqref{eq:advection}. We start by performing tests in 2D, for which the numerical results presented are obtained using two methods: 
	\begin{enumerate}
		\item polygonal ELLAM, obtained by approximating the cells and their trace-back regions with polygons (see, e.g. Figure \ref{trace-back-regions}, right), with mass conservation achieved approximately, with a relative error less than $10^{-4}$, by performing volume adjustments as in \cite{CDL18-GEM}. Here, the polygonal intersections are computed using a general polygon clipper (GPC) library, obtained from \url{http://www.cs.man.ac.uk/~toby/gpc/}.
		\item B-char ELLAM, as described in Section \ref{sec:B-char}. For this method, an $N$ has to be chosen to stop the iterations \eqref{eq:volAdj} before solving the optimisation problem \eqref{eq:opti_masscons}. Figure \ref{plot_N_massCons} shows the relative error on the mass balances against $N$ for a typical test case, and indicates that a reasonable choice is $N=10$ (reducing the errors to about 5\%). Further reduction does not bring much improvement.
	\end{enumerate} 
	For the first three test cases, we seek the concentration at time $T=8$, i.e. $c(\x,8)$. We also assume that the velocity field $\mathbf{u}$ is provided. These simple test cases aim at showing that the B-char ELLAM can achieve, with a much cheaper computational cost, numerical results that are essentially the same as those obtained from the polygonal ELLAM. Following these, we then use the B-char ELLAM for some 2D benchmark test cases, to show the robustness of the numerical scheme. Simple numerical tests (similar to the first three test cases in 2D) are then performed to show the applicability and efficiency of the B-char ELLAM in 3D. 
	
The relative errors will be measured in the $L^1$ and $L^2$ norm, by providing for $p=1,2$ the quantities
\[
E_p:= \dfrac{\norm{\Pi_\discC c(\x,T)-c(\x,T)}{p}}{\norm{ c(\x,T)}{p}},
\]
where $\norm{\cdot}{p}$ denote the norm in $L^p(\O)$.
	
	\subsection{Numerical tests in 2D}\label{sec:numTest2D}
	The test cases in 2D are performed over the domain $\O=(0,1)\times(0,1)$.

	\subsubsection{Test case 1}\label{2D:tc1}
	For the first test case, we consider a simple velocity field that simulates a translation along the $x$ axis, $\mathbf{u}=(\frac{1}{16},0)$. Although the no-flow boundary conditions are not satisfied, the final time $T=8$ is small enough so that no relevant characteristic traces outside the domain. The initial condition is set to be
	\begin{equation}\nonumber
	c(\x,0) =\begin{cases} & 1 \qquad \text{ if } \frac{1}{16}\leq x \leq \frac{5}{16}, \frac{1}{16}\leq y \leq \frac{5}{16} \\
	& 0 \qquad \text{ elsewhere }
	\end{cases}.
	\end{equation}
	Based on this initial condition and the given velocity, we expect the square block initially on the lower left corner of the domain to be transferred to the lower right corner of the domain, i.e. 
	\begin{equation}\nonumber
	c(\x,8) =\begin{cases} & 1 \qquad \text{ if } \frac{9}{16}\leq x \leq \frac{13}{16}, \frac{1}{16}\leq y \leq \frac{5}{16} \\
	& 0 \qquad \text{ elsewhere }
	\end{cases}.
	\end{equation}
	 
\begin{figure}[h]
		\centering
		\begin{tabular}{c@{\hspace*{2em}}c}
			\includegraphics[width=0.4\textwidth]{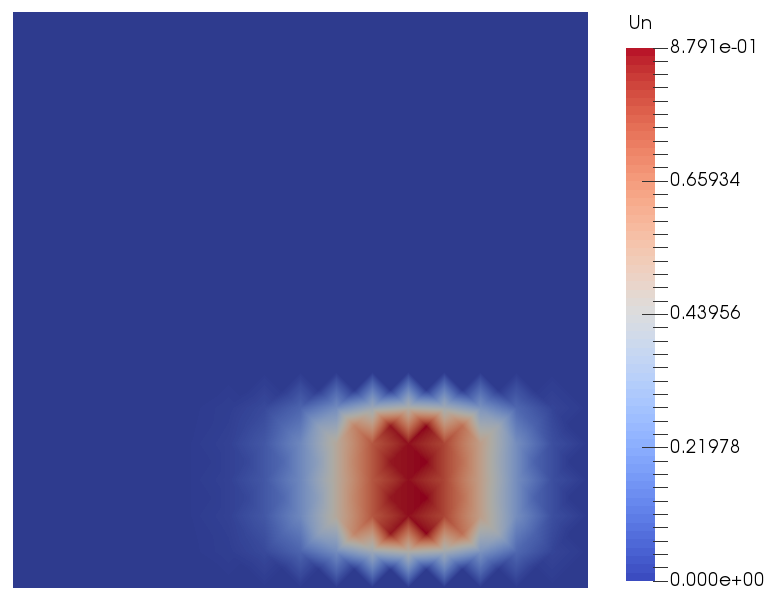} & 		 \includegraphics[width=0.4\textwidth]{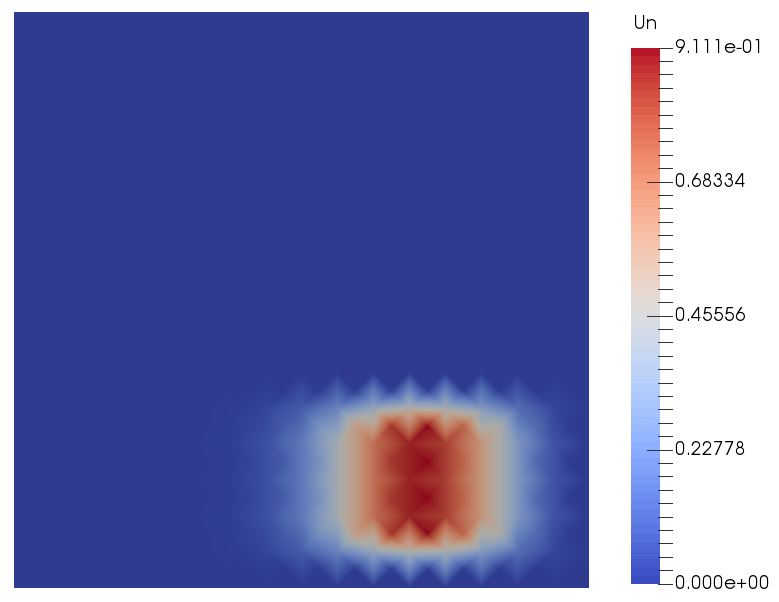}\\
		\end{tabular}
		\caption{ Concentration profiles obtained at final time $T=8$ with $\delta t = 0.8$ using an ELLAM scheme, test case 1 (left: polygonal; right: B-char).}
		\label{ELLAM_trans_dtpt8_mesh23}
	\end{figure}

\begin{figure}[h]
		\centering
		\begin{tabular}{c@{\hspace*{2em}}c}
			\includegraphics[width=0.4\textwidth]{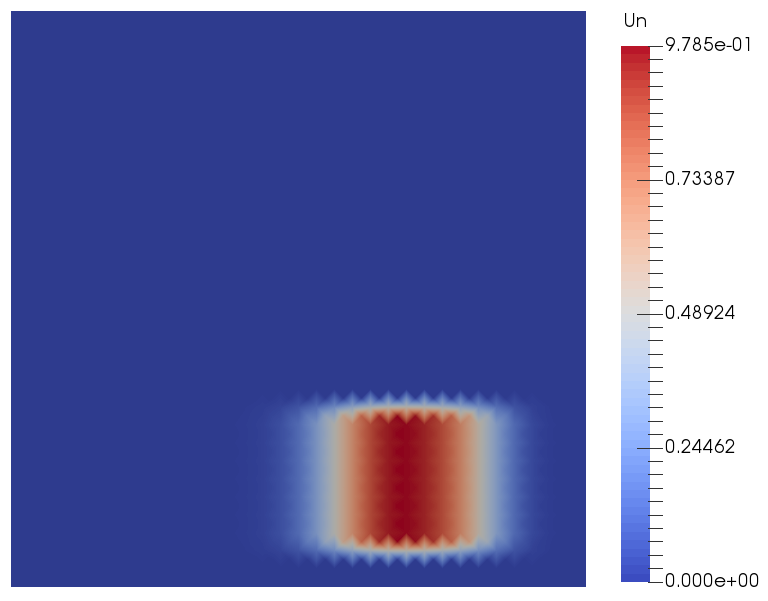} & 		 \includegraphics[width=0.4\textwidth]{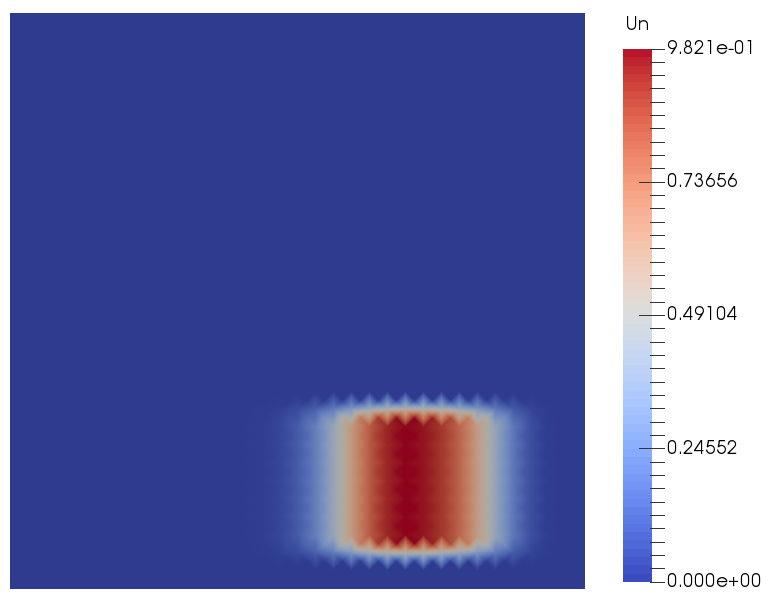}\\
		\end{tabular}
		\caption{ Concentration profiles obtained at final time $T=8$ with $\delta t = 0.4$ using an ELLAM scheme, test case 1 (left: polygonal; right: B-char).}
		\label{ELLAM_trans_dtpt4_mesh24}
	\end{figure}

\begin{figure}[h]
		\centering
		\begin{tabular}{c@{\hspace*{2em}}c}
			\includegraphics[width=0.4\textwidth]{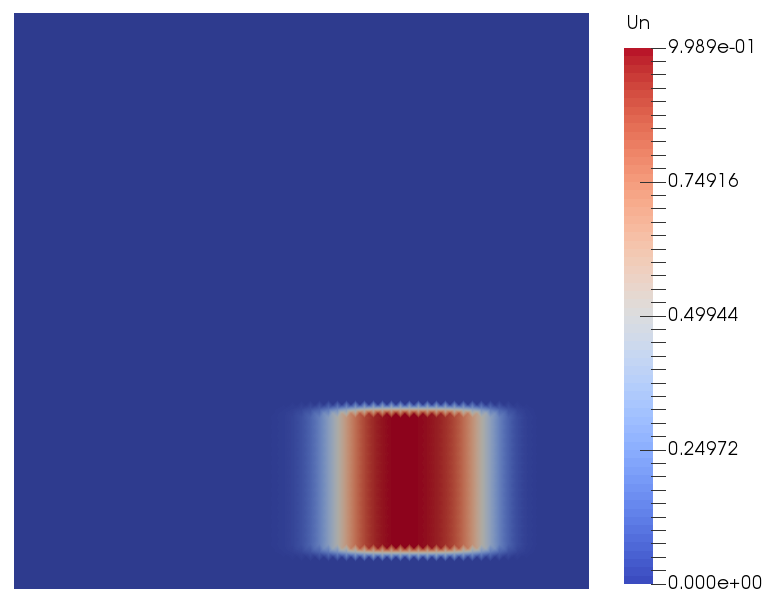} & 		 \includegraphics[width=0.4\textwidth]{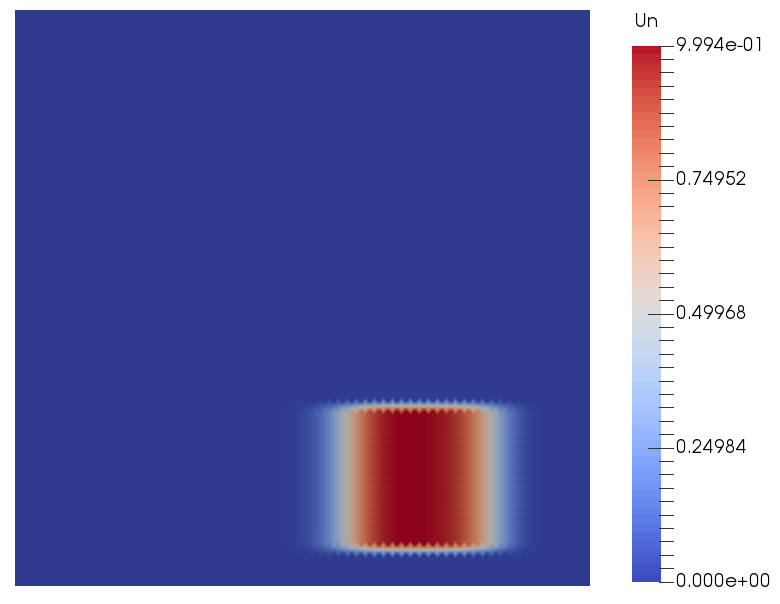}\\
		\end{tabular}
		\caption{ Concentration profiles obtained at final time $T=8$ with $\delta t = 0.2$ using an ELLAM scheme, test case 1 (left: polygonal; right: B-char).}
		\label{ELLAM_trans_dtpt2_mesh25}
	\end{figure}

	We now compute and compare the approximate solutions using the polygonal ELLAM scheme and the B-char ELLAM. For the B-char ELLAM, 4 balls are used to approximate each cell. This will be performed starting on a $16\times16$ grid with a time step of $\delta t= 0.8$, and refined for 2 levels in space and time, leading to a test on a $64\times64$ grid with a time step of $\delta t =0.2$. Upon looking at Figures \ref{ELLAM_trans_dtpt8_mesh23} to \ref{ELLAM_trans_dtpt2_mesh25}, we see that the concentration profiles obtained from the polygonal ELLAM and the B-char ELLAM are quite similar, with the B-char ELLAM producing maximum concentrations which are slightly closer to 1, as compared to the polygonal ELLAM.  
	
	Now, we compare these methods in more detail by looking at Tables \ref{tab:transErr-std} and \ref{tab:transErr-ball}. As can be seen, the polygonal ELLAM and the B-char ELLAM produce results that are quite close to one another. Moreover, upon measuring the CPU runtime (in seconds) for one time step for the total process of tracking, computing intersections, and performing volume adjustments, we see that the B-char ELLAM gets to perform the simulations much faster compared to the polygonal ELLAM. This is mainly due to the fact that ball intersections are much cheaper to compute as compared to polygonal intersections. 
	
	\begin{table}[h]
		\begin{center}
			\begin{tabular}{|c|c|c|c|c|}
				\hline
				Mesh & $\delta t$ & CPU time & $E_1$ &  $E_2$  \\
				 & & (one time step) &  &  \\
				\hline
				$16\times16$ & 0.8 & 0.5175 & 4.8271e-01 & 3.7277e-01 \\
				\hline
				$32\times32$  & 0.4 & 6.4640 & 3.4911e-01 & 3.1673e-01\\
				\hline
				$64\times64$ & 0.2 & 97.3994 & 2.4956e-01 & 2.6898e-01\\
				\hline
			\end{tabular} 
		\end{center}
		\caption{CPU runtime and errors in the concentration profiles, test case 1, polygonal ELLAM, final time $T=8$.}\label{tab:transErr-std}
	\end{table} 
		
	\begin{table}[h]
		\begin{center}
			\begin{tabular}{|c|c|c|c|c|}
				\hline
				Mesh & $\delta t$ & CPU time & $E_1$ &  $E_2$  \\
				& & (one time step) &  &  \\
				\hline
				$16\times16$ & 0.8 & 0.1141 & 4.7637e-01 & 3.8273e-01 \\
				\hline
				$32\times32$  & 0.4 & 0.4321 & 3.4889e-01 & 3.3183e-01\\
				\hline
				$64\times64$ & 0.2 & 3.5188 & 2.5558e-01 & 2.9220e-01\\
				\hline
			\end{tabular} 
		\end{center}
		\caption{CPU runtime and errors in the concentration profiles, test case 1, B-char, final time $T=8$.}\label{tab:transErr-ball}
	\end{table}

\begin{remark}[CPU runtime] The CPU times are only used as an indication to show the advantage of the B-char method over the polygonal ELLAM. The codes may not be fully optimised, but are implemented in a similar manner for both methods, by taking advantage of the vectorial capacities of MATLAB.
\end{remark}

 \subsubsection{Test case 2}\label{2D:tc2}
\begin{figure}[h]
		\centering
		\includegraphics[width=0.8\textwidth]{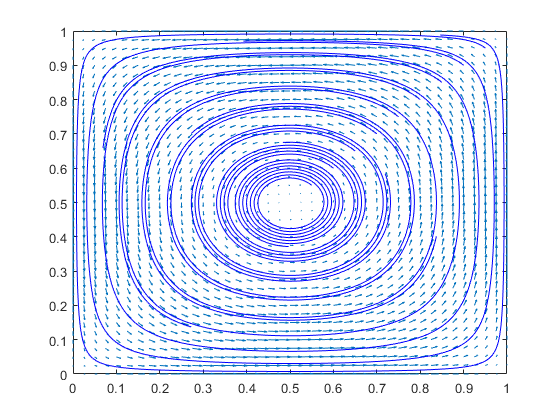} \\
		\caption{Streamlines of the velocity field $\mathbf{u}=((1-2y)(x-x^2),-(1-2x)(y-y^2))$.}
		\label{exactVel}
	\end{figure}

	The second test case considers a velocity field $\mathbf{u}=((1-2y)(x-x^2),-(1-2x)(y-y^2))$. Here, $\mathbf{u}$ is a divergence-free velocity field which simulates a rotation with some stretching, and the centre of this rotation is located at $(0.5,0.5)$ (see Figure \ref{exactVel}). 
	
\begin{figure}[h]
		\centering
		\begin{tabular}{c@{\hspace*{2em}}c}
			\includegraphics[width=0.4\textwidth]{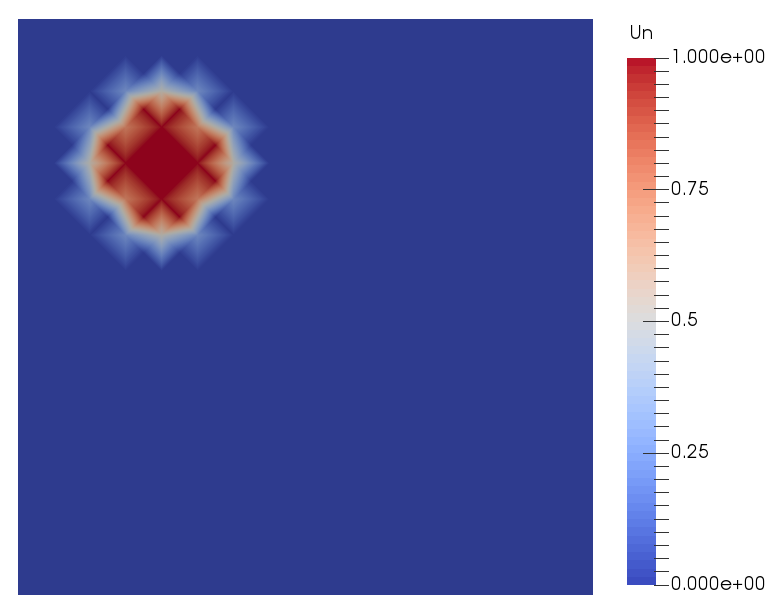} & 		\includegraphics[width=0.4\textwidth]{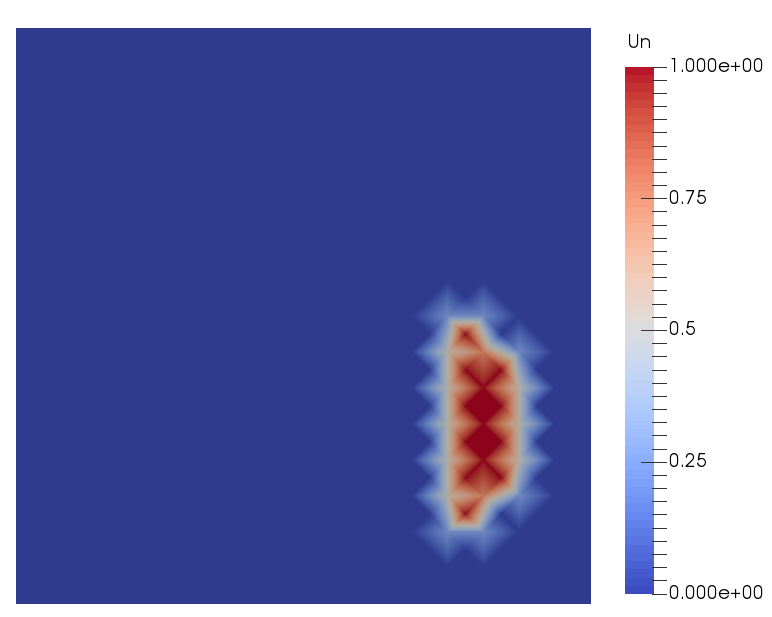}\\
		\end{tabular}
		\caption{ $c(\x,t)$ for test case 2 (left: initial condition at $t=0$; right: benchmark solution profile at final time $T=8$).}
		\label{exactSol}
	\end{figure}

\begin{figure}[h]
		\centering
		\begin{tabular}{c@{\hspace*{2em}}c}
			\includegraphics[width=0.4\textwidth]{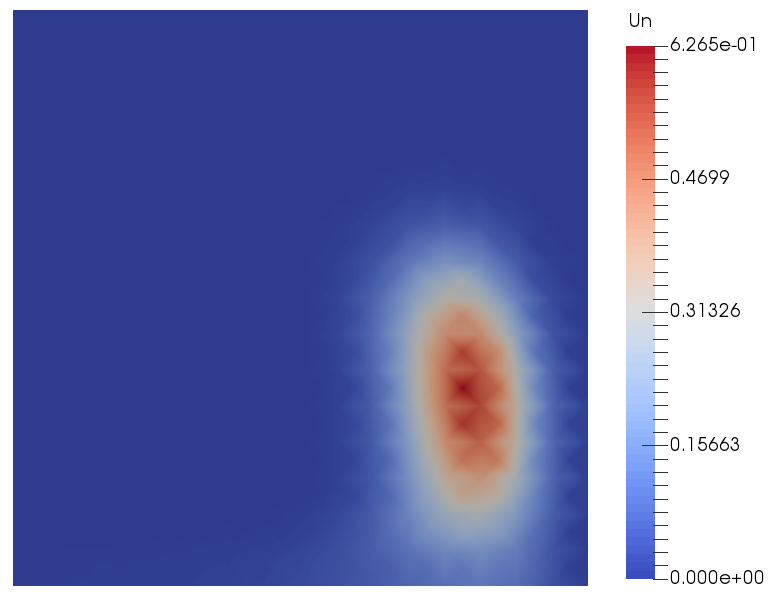} & 		 \includegraphics[width=0.4\textwidth]{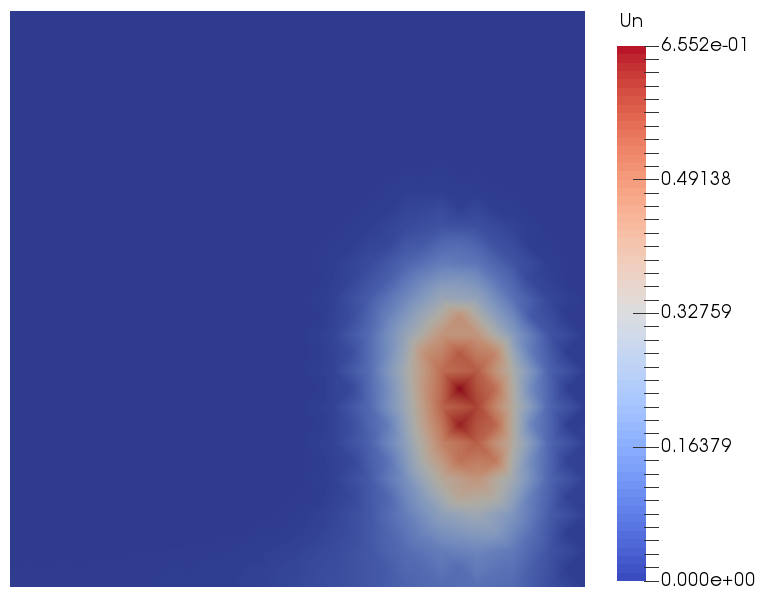}\\
		\end{tabular}
		\caption{ Concentration profiles obtained at final time $T=8$ with $\delta t = 0.8$ using an ELLAM scheme, test case 2 (left: polygonal; right: B-char).}
		\label{ELLAM_rot_dtpt8_mesh23}
	\end{figure}

\begin{figure}[h]
		\centering
		\begin{tabular}{c@{\hspace*{2em}}c}
			\includegraphics[width=0.4\textwidth]{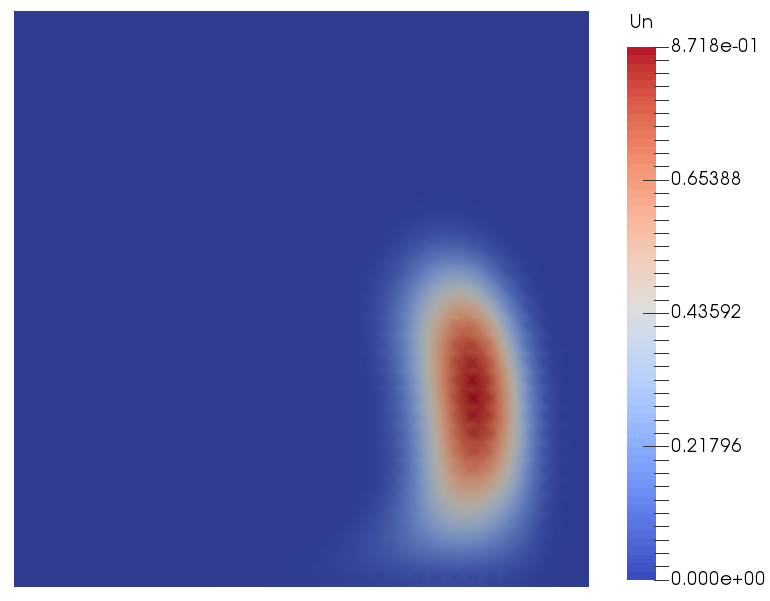} & 		 \includegraphics[width=0.4\textwidth]{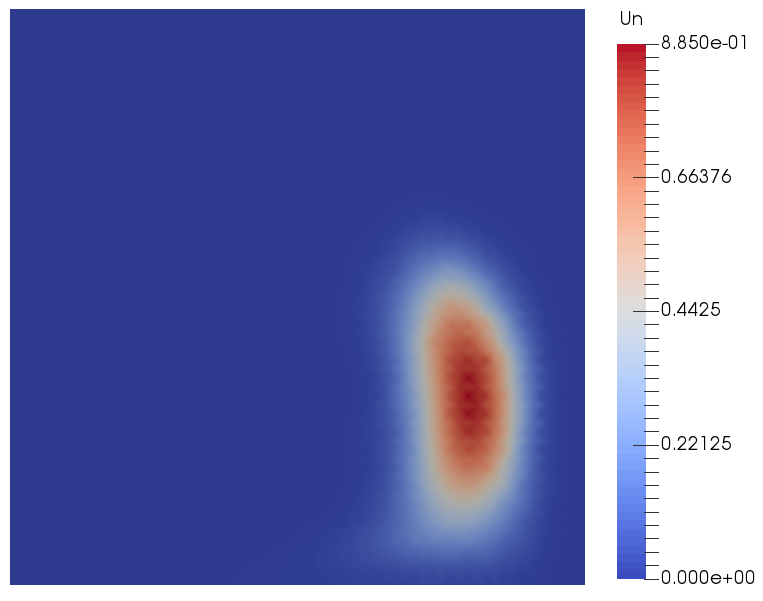}\\
		\end{tabular}
		\caption{ Concentration profiles obtained at final time $T=8$ with $\delta t = 0.4$ using an ELLAM scheme, test case 2 (left: polygonal; right: B-char).}
		\label{ELLAM_rot_dtpt4_mesh24}
	\end{figure}

\begin{figure}[h]
		\centering
		\begin{tabular}{c@{\hspace*{2em}}c}
			\includegraphics[width=0.4\textwidth]{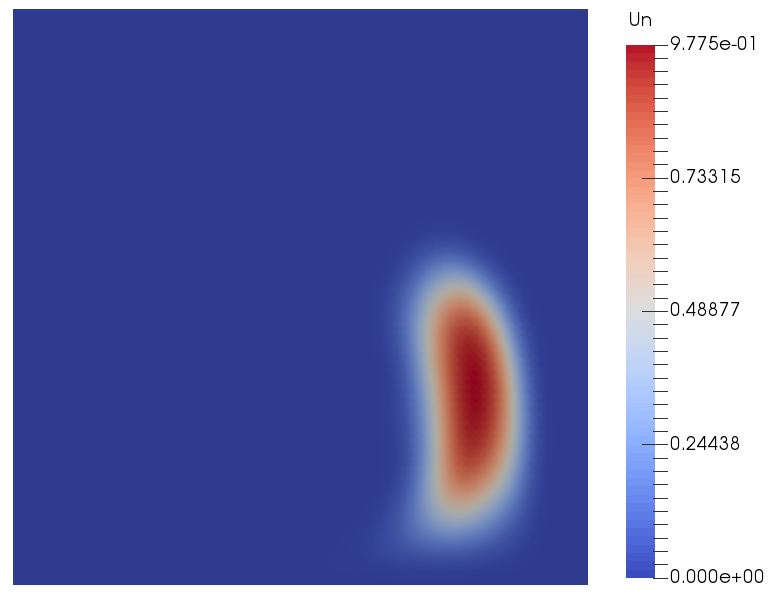} & 		 \includegraphics[width=0.4\textwidth]{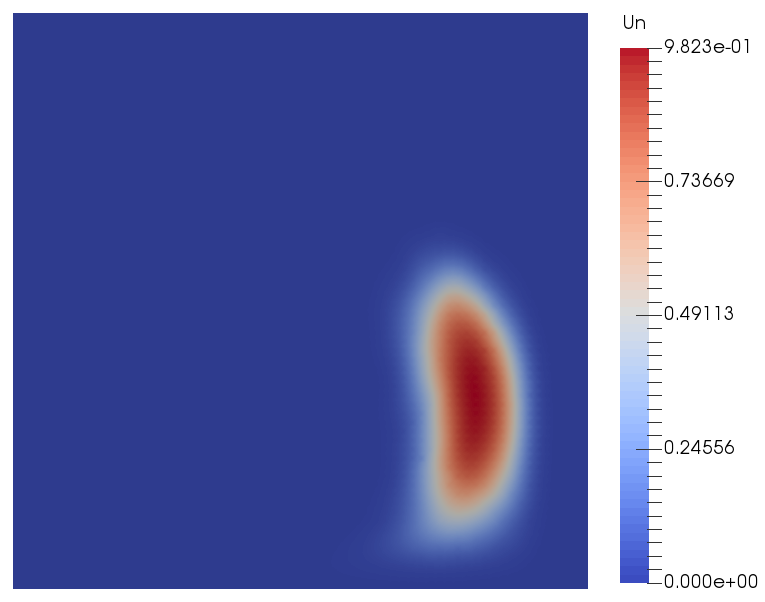}\\
		\end{tabular}
		\caption{ Concentration profiles obtained at final time $T=8$ with $\delta t = 0.2$ using an ELLAM scheme, test case 2 (left: polygonal; right: B-char).}
		\label{ELLAM_rot_dtpt2_mesh25}
	\end{figure}

	The initial condition is set to be
	\begin{equation}\nonumber
	c(\x,0) =\begin{cases} & 1 \qquad \text{ if } (x-\frac{1}{4})^2+(y-\frac{3}{4})^2 < \frac{1}{64} \\
	& 0 \qquad \text{ elsewhere }
	\end{cases}.
	\end{equation} 
	Essentially, this assumes that we have a substance near the top-left corner of our domain (see Figure \ref{exactSol}, left), being rotated, and somehow stretched for $t=8$ time units. Unlike the first test case, an exact solution is not available. Hence, we compare our results with a benchmark solution, obtained by solving \eqref{charac} using an Euler method over a very fine grid (to be particular, 2 levels of refinement compared to the mesh being considered), with a very small time step $\delta t = 0.001$. This is then projected onto the mesh being considered$-$in the case of Figure \ref{exactSol}, right, a mesh consisting of $16\times16$ squares.
	As with the first test case, we start by comparing the concentration profiles obtained by solving \eqref{eq:advection} using the polygonal ELLAM and the B-char ELLAM, with 4 balls being used to approximate each cell for the B-char ELLAM. Based on Figures \ref{ELLAM_rot_dtpt8_mesh23} to \ref{ELLAM_rot_dtpt2_mesh25}, we see that the concentration profile obtained from the B-char ELLAM is very similar to those obtained from the polygonal ELLAM. Also, as with the first test case, the maximum concentrations for the B-char ELLAM are closer to 1, compared to the polygonal ELLAM. 
	
	Upon performing a more rigorous comparison by looking at Tables \ref{tab:rotErr-std} and \ref{tab:rotErr-ball}, we note that the errors in both the $L^1$ and $L^2$ norm for both methods are quite close to each other. Also, the B-char ELLAM performs much faster than the polygonal ELLAM. 
	
	\begin{table}[h]
		\begin{center}
			\begin{tabular}{|c|c|c|c|c|}
				\hline
				Mesh & $\delta t$ & CPU time & $E_1$ &  $E_2$   \\
				& & (one time step) &  &  \\
				\hline
				$16\times16$ & 0.8 & 2.7441 & 7.3431e-01 & 5.1027e-01 \\
				\hline
				$32\times32$  & 0.4 & 43.0472 & 6.3375e-01 & 4.2258e-01\\
				\hline
				$64\times64$ & 0.2 & 700.9942 & 4.9580e-01 & 3.6537e-01\\
				\hline
			\end{tabular} 
		\end{center}
		\caption{CPU runtime and errors in the concentration profiles, test case 2, polygonal ELLAM, final time $T=8$.}\label{tab:rotErr-std}
	\end{table} 
	
	\begin{table}[h]
		\begin{center}
			\begin{tabular}{|c|c|c|c|c|}
				\hline
				Mesh & $\delta t$ & CPU time & $E_1$ &  $E_2$  \\
				& & (one time step) &  &  \\
				\hline
				$16\times16$ & 0.8 & 0.1865 & 7.3138e-01 & 5.0673e-01 \\
				\hline
				$32\times32$  & 0.4 & 1.3095 & 6.1391e-01 & 4.1428e-01\\
				\hline
				$64\times64$ & 0.2 & 14.5061 & 4.7916e-01 & 3.5931e-01\\
				\hline
			\end{tabular} 
		\end{center}
		\caption{CPU runtime and errors in the concentration profiles, test case 2, B-char, final time $T=8$.}\label{tab:rotErr-ball}
	\end{table} 

 \subsubsection{Test case 3}\label{2D:tc3}
\begin{figure}[h]
	\centering
	\begin{tabular}{c@{\hspace*{2em}}c}
		\includegraphics[width=0.4\textwidth]{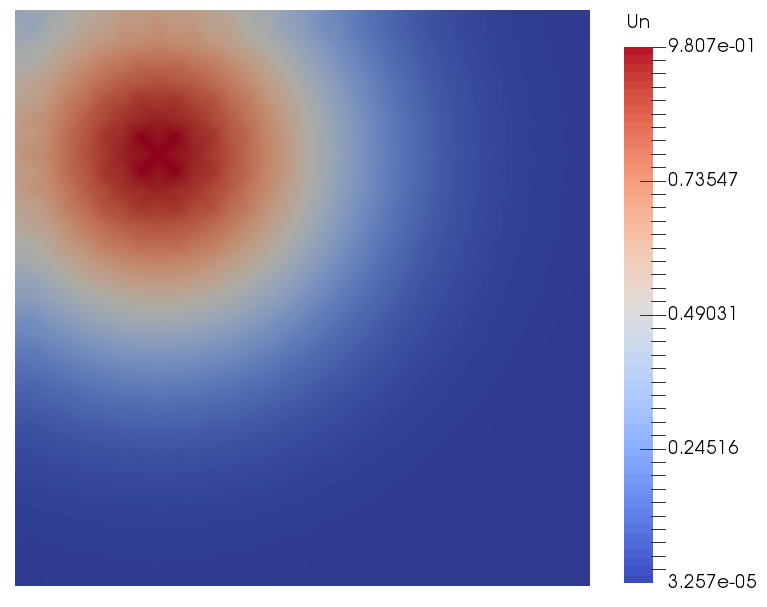} & 		\includegraphics[width=0.4\textwidth]{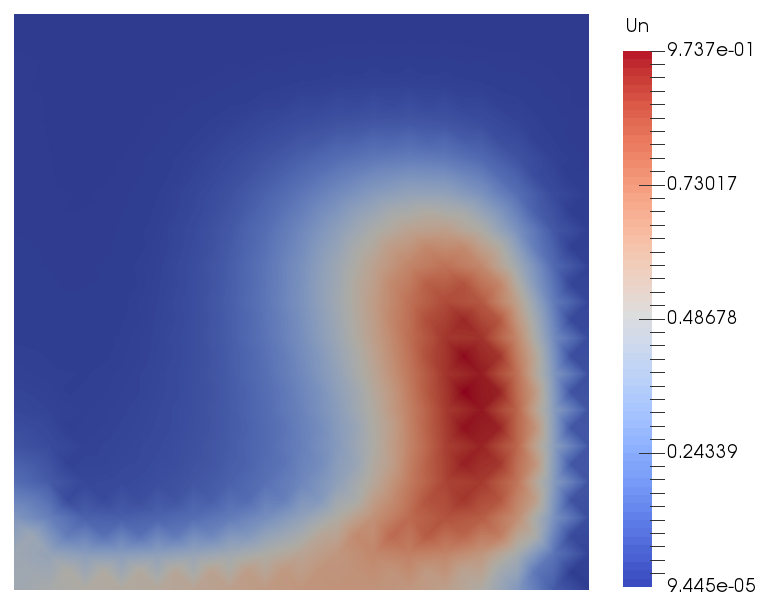}\\
	\end{tabular}
	\caption{$c(\x,t)$ for test case 3 (left: initial condition at $t=0$; right: benchmark solution profile at final time $T=8$).}
	\label{exactSol-cont}
\end{figure}

 Now, we present a test case using the same velocity field as the second test, but now with a smooth initial condition, given by $c(\x,0) = \exp(-10((x-\frac{1}{4})^2+(y-\frac{3}{4})^2))$. Similar to the second test case, this assumes that we have majority of our substance near the top-left corner of our domain (see Figure \ref{exactSol-cont}, left). The benchmark solution is also obtained in the same way as the second test case, and is then projected onto a mesh consisting of $16\times16$ squares (Figure \ref{exactSol-cont}, right).
	
	\begin{figure}[h]
		\centering
		\begin{tabular}{c@{\hspace*{2em}}c}
			\includegraphics[width=0.4\textwidth]{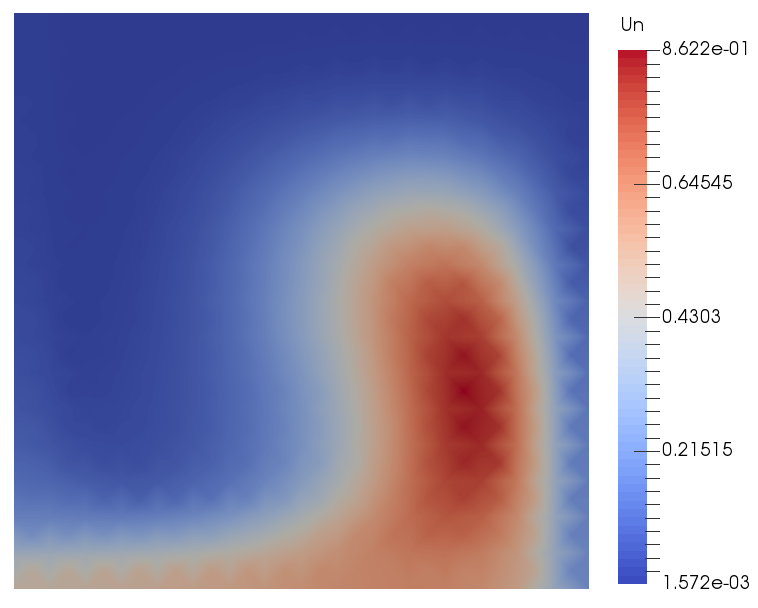} & 		 \includegraphics[width=0.4\textwidth]{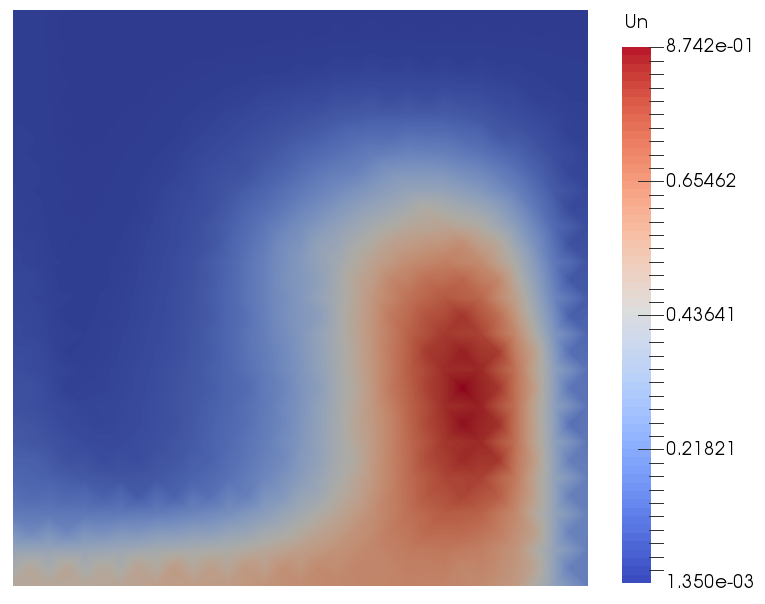}\\
		\end{tabular}
		\caption{ Concentration profiles obtained at final time $T=8$ with $\delta t = 0.8$ using an ELLAM scheme, test case 3 (left: polygonal; right: B-char).}
		\label{ELLAM_cont_rot_dtpt8_mesh23}
	\end{figure}
	
	\begin{figure}[h]
		\centering
		\begin{tabular}{c@{\hspace*{2em}}c}
			\includegraphics[width=0.4\textwidth]{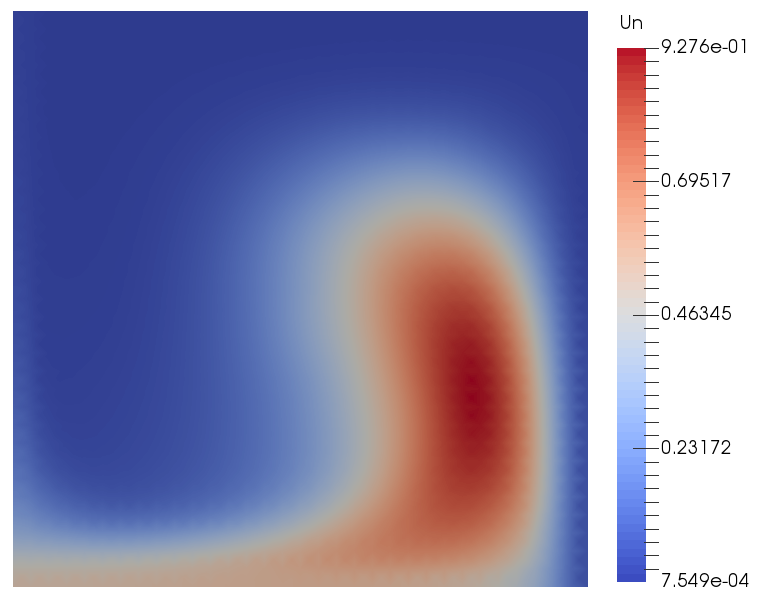} & 		 \includegraphics[width=0.4\textwidth]{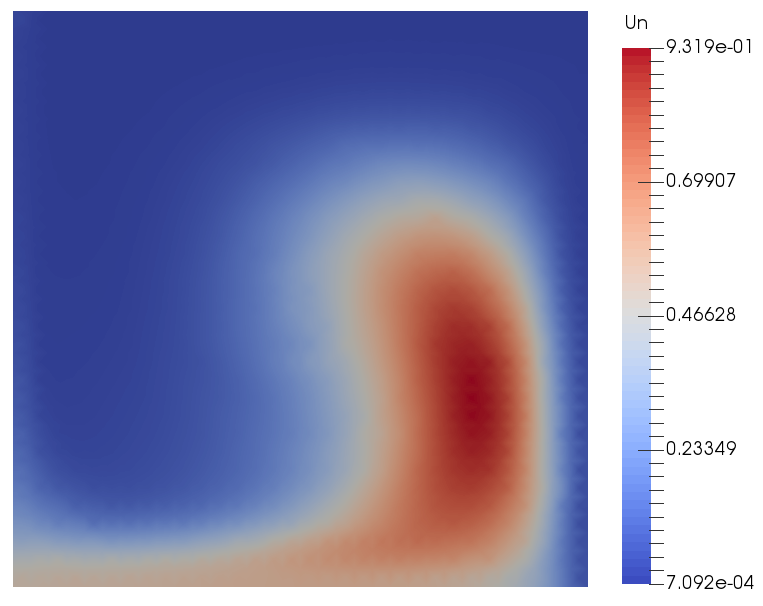}\\
		\end{tabular}
		\caption{ Concentration profiles obtained at final time $T=8$ with $\delta t = 0.4$ using an ELLAM scheme, test case 3 (left: polygonal; right: B-char).}
		\label{ELLAM_cont_rot_dtpt4_mesh24}
	\end{figure}
	
	\begin{figure}[h]
		\centering
		\begin{tabular}{c@{\hspace*{2em}}c}
			\includegraphics[width=0.4\textwidth]{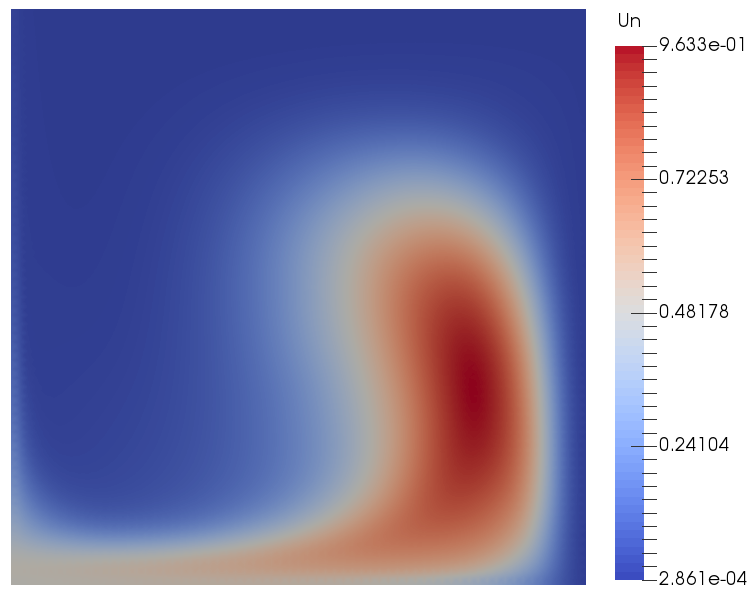} & 		 \includegraphics[width=0.4\textwidth]{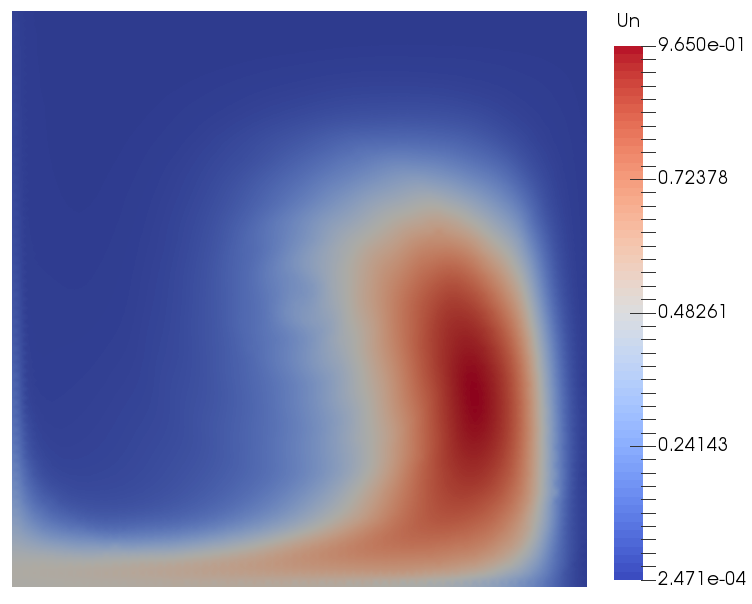}\\
		\end{tabular}
		\caption{Concentration profiles obtained at final time $T=8$ with $\delta t = 0.2$ using an ELLAM scheme, test case 3 (left: polygonal; right: B-char).}
		\label{ELLAM_cont_rot_dtpt2_mesh25}
	\end{figure}
	
	From Figures \ref{ELLAM_cont_rot_dtpt8_mesh23} to \ref{ELLAM_cont_rot_dtpt2_mesh25} and Tables \ref{tab:rotErr-cont-std} and \ref{tab:rotErr-cont-ball}, a similar observation can be made as with the first two test cases: the concentration profiles obtained from the B-char ELLAM are very close to those obtained from the polygonal ELLAM, and the B-char ELLAM also performs much faster than the polygonal ELLAM. It is also notable that due to the continuous initial condition, the errors for this test case are smaller than those obtained from the second test case (by more than a factor of 5). 
	
	\begin{table}[h]
		\begin{center}
			\begin{tabular}{|c|c|c|c|c|}
				\hline
				Mesh & $\delta t$ & CPU time & $E_1$ &  $E_2$  \\
				& & (one time step) &  &  \\
				\hline
				$16\times16$ & 0.8 & 2.7276 & 1.4141e-01 & 1.4244e-01 \\
				\hline
				$32\times32$  & 0.4 & 42.6126 & 9.3292e-02 & 9.9508e-02\\
				\hline
				$64\times64$ & 0.2 & 724.9937 & 5.6584e-02 & 6.7150e-02\\
				\hline
			\end{tabular} 
		\end{center}
		\caption{CPU runtime and errors in the concentration profiles, test case 3, polygonal ELLAM, final time $T=8$.}\label{tab:rotErr-cont-std}
	\end{table} 
	
	\begin{table}[h]
		\begin{center}
			\begin{tabular}{|c|c|c|c|c|}
				\hline
				Mesh & $\delta t$ & CPU time & $E_1$ &  $E_2$  \\
				& & (one time step) &  &  \\
				\hline
				$16\times16$ & 0.8 & 0.2146 & 1.4961e-01 & 1.5055e-01 \\
				\hline
				$32\times32$  & 0.4 & 1.3015 & 9.1979e-02 & 9.8428e-02\\
				\hline
				$64\times64$ & 0.2 & 16.0974 & 5.6735e-02 & 6.7733e-02\\
				\hline
			\end{tabular} 
		\end{center}
		\caption{CPU runtime and errors in the concentration profiles, test case 3, B-char, final time $T=8$.}\label{tab:rotErr-cont-ball}
	\end{table} 
	
	In general, we see that for both the polygonal ELLAM scheme and B-char ELLAM, numerical diffusion is most prominent on the coarse $16\times16$ grid. This numerical diffusion becomes less prominent as the grid and the time step are refined. Moreover, for all test cases, using the B-char ELLAM achieves the same level of accuracy as the polygonal ELLAM, while at the same time, reducing the required computational cost (up to 40 times faster on the finest mesh in our tests). 
	
\subsubsection{Solid body rotation test case} \label{subsubsec:Solidbodyrot}
\begin{figure}[h]
	\centering
	\begin{tabular}{c@{\hspace*{2em}}c}
		\includegraphics[width=0.4\textwidth]{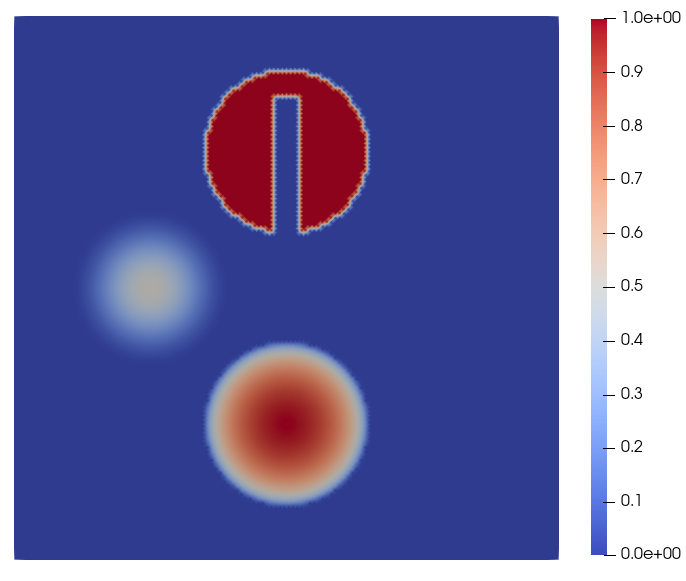} & 		 \includegraphics[width=0.4\textwidth]{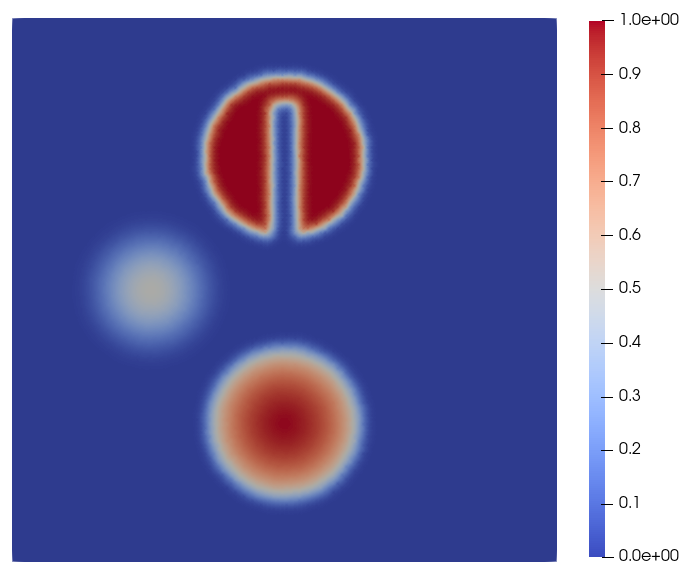}\\
	\end{tabular}
	\caption{ Concentration profile for the solid body rotation test case (left: initial condition; right: numerical solution at $T=2\pi$).}
	\label{fig:solidBodyRot2D}
\end{figure}

	Now that we have established that the B-char ELLAM achieves the same level of accuracy as the polygonal ELLAM with a cheaper computational cost, we perform a numerical simulation on a benchmark test case, the solid body rotation problem proposed in \cite{L96-benchmarkTestcase}. Instead of the domain $\O= (0,1)\times (0.5,1.5)$ considered therein, we retain our domain at $\O=(0,1) \times (0,1)$, and modify the initial conditions accordingly, as done in \cite{BRP14}. To be specific, the initial condition is made up of three components: a bump described by 
\begin{equation}\nonumber
\begin{aligned}
q(\x) &= 0.25 (1+\cos(\pi r(x,y))),\\
r(\x) &= \frac{1}{0.15}\min(\sqrt{(x-0.25)^2 + (y-0.5)^2},0.15),
\end{aligned}
\end{equation}
together with a cone and a slotted cylinder of radius $0.15$ and height $1$, centered at $(0.5,0.25)$ and $(0.5,0.75)$, respectively. The slot in the cylinder is created by removing the region $[0.475,0.525]\times[0.6,0.85]$ (see Figure \ref{fig:solidBodyRot2D}, left). For this test case, the velocity field is given by $\mathbf{u} = ((0.5-y),(x-0.5))$; we note that it does not satisfy the no-flow boundary conditions $\mathbf{u}\cdot\mathbf{n}=0$ on the boundary of the domain, but this has no impact on the B-char method since none of the balls that intersect non-zero values of the solution are tracked near the boundary. This velocity field simulates a rotation, and completes a revolution at time $T=2\pi$. Hence, we expect the solution profile at time $T=2\pi$ to be identical to the initial condition. Numerical simulations are performed using the B-char ELLAM with 4 balls in each cell. As with the test cases provided in \cite{BRP14,K06}, we take a mesh with $128\times 128$ cells, but to take advantage of the ELLAM, we use a much larger time step of $\delta t = 2\pi/10$, compared to $\delta t = 2\pi/810$ or smaller in the literature.


 As can be seen in Figure \ref{fig:solidBodyRot2D}, the numerical solution obtained after one cycle of revolution is almost the same as the initial condition. In particular, the computed relative errors are $E_1 = 1.3630$e-$01$ and $E_2 = 2.1178$e-$01$. The absolute errors obtained in \cite{BRP14} (using a different scheme) in $L^1$ and $L^2$ norms are $1.7$e-$02$ and $7.5$e-$02$, respectively. These absolute errors correspond to relative errors of $1.4560$e-$01$ and $2.4020$e-$01$, respectively. These show that, despite its cheap computational cost, the B-char ELLAM method produces solutions with an accuracy that is comparable to those obtained in the literature. Moreover, the numerical solution preserves the slot in the cylinder, which indicates that the numerical diffusion introduced by the B-char ELLAM is minimal. 

\subsubsection{Deformational flow test} 
\begin{figure}[h]
		\centering
		\begin{tabular}{c@{\hspace*{2em}}c}
			\includegraphics[width=0.4\textwidth]{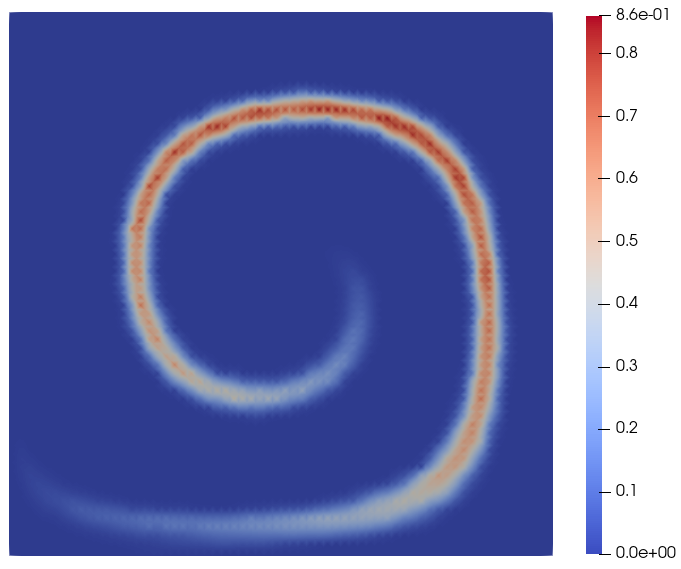} & 		 \includegraphics[width=0.4\textwidth]{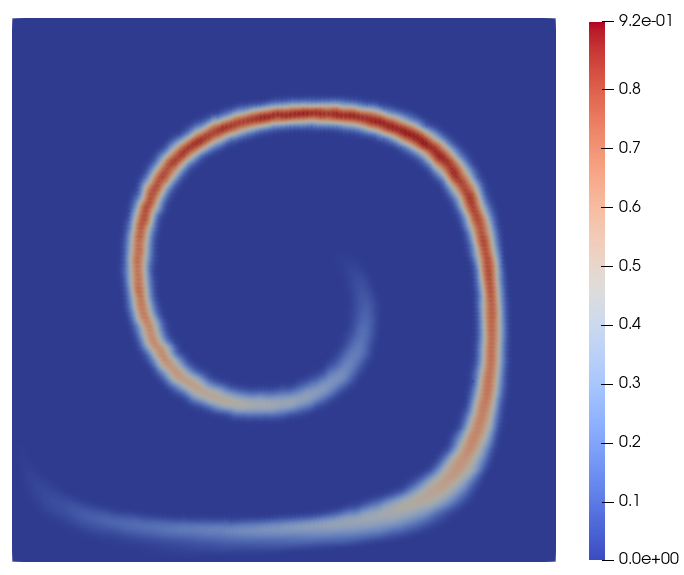}\\
		\end{tabular}
		\caption{ Concentration profile for the deformational flow test at halftime $T=2.5$  (left: $64 \times 64$ cells; right: $128 \times 128$ cells).}
		\label{fig:deformationalFlowT2}
	\end{figure}

\begin{figure}[h]
		\centering
		\begin{tabular}{c@{\hspace*{2em}}c}
			\includegraphics[width=0.4\textwidth]{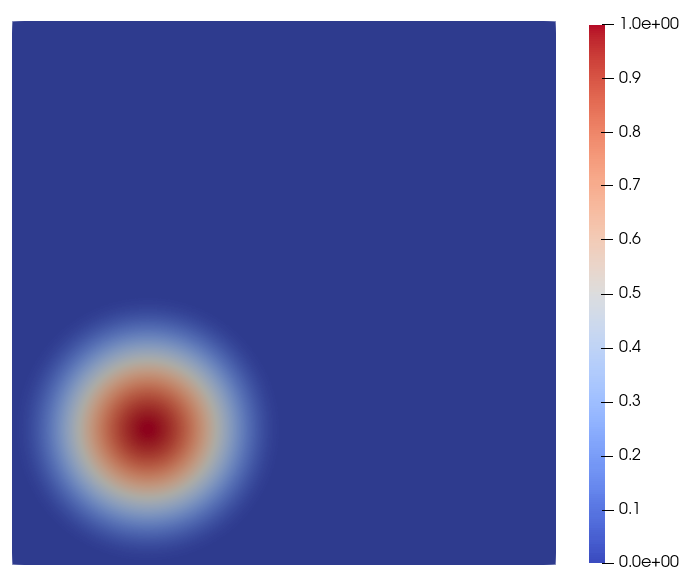} & 		 \includegraphics[width=0.4\textwidth]{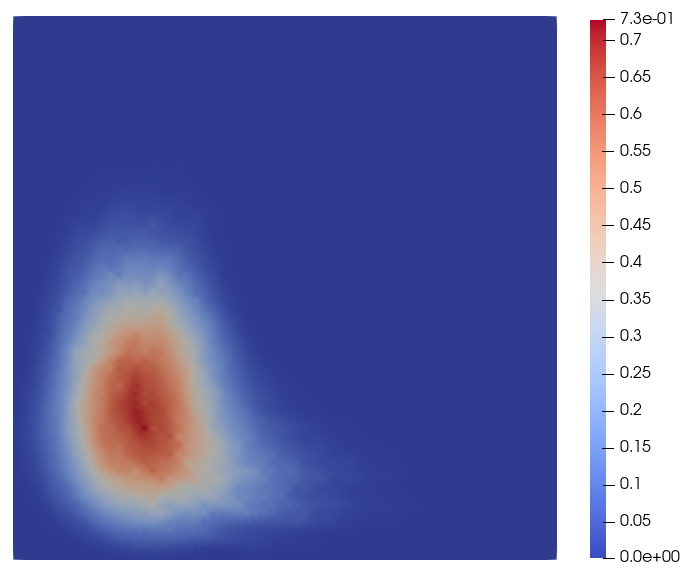}\\
		\end{tabular}
		\caption{ Concentration profile for the deformational flow test, $64 \times 64$ cells (left: initial condition; right: numerical solution at $T=5$).}
		\label{fig:deformationalFlow25}
	\end{figure}

\begin{figure}[h]
	\centering
	\begin{tabular}{c@{\hspace*{2em}}c}
		\includegraphics[width=0.4\textwidth]{initCond_deformationaltest_2D} & 		 \includegraphics[width=0.4\textwidth]{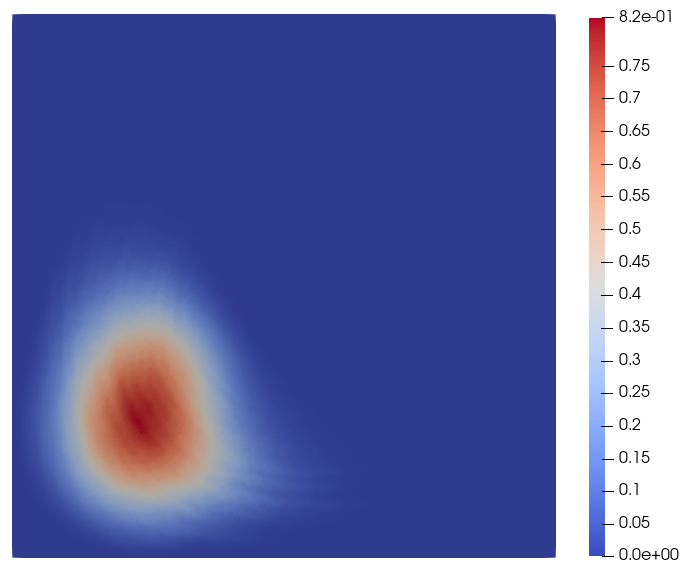}\\
	\end{tabular}
	\caption{ Concentration profile for the deformational flow test, $128 \times 128$ cells (left: initial condition; right: numerical solution at $T=5$).}
	\label{fig:deformationalFlow26}
\end{figure}

Finally, we present a test case for which the velocity field gives a very strong deformation. In  \cite{LT12-deformationtest}, such a velocity field was given in terms of spherical coordinates, but \cite{NL10-deformationtest} provided an analogue in Cartesian coordinates:
\[
\mathbf{u} = (\sin^2(\pi x)\sin(2\pi y)\cos(\pi t/T),-\sin^2(\pi y)\sin(2\pi x)\cos(\pi t/T)).
\]
This velocity field simulates a deformational flow, which reaches zero and then changes direction at half time $T/2$, resulting in a reversal of the flow field. We expect the solution profile to be the most highly deformed at $T/2$, and to be identical to the initial condition at time $T$. Following the test proposed in \cite[Section 5.7.4]{D98-NumericalFluidDynamics}, we set $T=5$ and take our initial condition to be a cosine bell function:
	\begin{equation}\nonumber
	\begin{aligned}
	c(\x,0) &= 0.5 (1+\cos(\pi r(\x))),\\
	r(\x) &= \min(4\sqrt{(x-0.25)^2 + (y-0.5)^2},1).
	\end{aligned}
	\end{equation} We perform numerical tests on the following meshes: a mesh with $64\times 64$ cells and a refinement which consists of $128\times 128$ cells. In order to take advantage of ELLAM, we take larger time steps of $\delta t = 0.5$ and $\delta t = 0.25$, respectively, compared to $\delta t = 0.01$ in the literature. As expected, in our numerical simulations, we obtain the most deformed solution profile at time $T=2.5$ (see Figure \ref{fig:deformationalFlowT2}). As noticed in Figures \ref{fig:deformationalFlow25} and \ref{fig:deformationalFlow26}, right, the strong distortion at time $T=2.5$ entails a distortion of the numerical solution at final time $T=5$. This distortion is more prominent for the coarse mesh. Aside from this, we notice the presence of some numerical diffusion, since the maximum amplitude of the solution has decreased, by around $27\%$ for the coarse mesh, and $18\%$ for the refined mesh. These issues were also encountered in \cite{D98-NumericalFluidDynamics} (with similar levels of distortion and diffusion), and can be resolved by further refinement of the mesh. 
	
	\subsection{Numerical tests in 3D}

	In 3D, we perform the numerical tests over the domain $\O=(0,1)\times(0,1)\times(0,1)$. These are only performed using the B-char ELLAM, by approximating each cell with 8 balls. The polygonal ELLAM was not used for the numerical tests here as an easy/efficient implementation of intersecting a convex and non-convex polyhedron is not readily available. The numerical tests are performed over a mesh with $16\times16\times16$ cubes, with a time step $\delta t = 0.8$. 

As with the 2D case, the first test case involves a translation about the $x$ axis by considering the velocity field $\mathbf{u}=(\frac{1}{16},0,0)$. The concentration is initially set to be at the cubic block on the lower left corner of the cube, given by
	\[
	c(\x,0) = \begin{cases} 1 \quad \text{ if }  \frac{1}{16}\leq x \leq \frac{5}{16},\frac{1}{16}\leq y \leq \frac{5}{16}, \frac{1}{16}\leq z \leq \frac{5}{16} \\
	0 \quad \text{ elsewhere }
	\end{cases}.
	\]

	For the second test case, we consider a translation along the $z$ axis, accompanied by some rotation and stretching along the $x$ and $y$ directions by considering the velocity field $\mathbf{u}=((1-2y)(x-x^2),-(1-2x)(y-y^2),\frac{1}{16})$. Initially, the concentration is located on a cylinder described by 
	\[
	c(\x,0) =\begin{cases}  1 \quad \text{ if } (x-\frac{1}{4})^2+(y-\frac{3}{4})^2 < \frac{1}{64}, \frac{1}{16}\leq z \leq\frac{5}{16} \\
	 0 \quad \text{ elsewhere }
	\end{cases}.
	\]
	
	Finally, the third test case also uses the velocity field $\mathbf{u}=((1-2y)(x-x^2),-(1-2x)(y-y^2),\frac{1}{16})$, with a continuous initial condition 
	\[
 c(\x,0) = \exp\left(-10\bigg(\big(x-\frac{1}{4}\big)^2+\big(y-\frac{3}{4}\big)^2+\big(z-\frac{3}{16}\big)^2\bigg)\right),
	\]
for which majority of the substance is concentrated around a neighborhood of the point $(\frac{1}{4},\frac{3}{4},\frac{3}{16})$. 
	
	We note that for all test cases, the no-flow boundary conditions are not satisfied, but the final time $T=8$ is such that none of the relevant characteristics trace outside the domain. For the first test case, an exact solution is available, given by
	\[
	c(\x,8) = \begin{cases} 1 \quad \mathrm{ if } \quad  \frac{9}{16}\leq x \leq \frac{13}{16},\frac{1}{16}\leq y \leq \frac{5}{16}, \frac{1}{16}\leq z \leq \frac{5}{16} \\
	0 \quad \mathrm{ elsewhere }
	\end{cases}.
	\] 
For the second and third test cases, benchmark solutions are computed using a similar method described in the second test case in Section \ref{sec:numTest2D}.
	\begin{table}[h]
		\begin{center}
			\begin{tabular}{|c|c|c|c|c|}
				\hline
				Test case & $\delta t$ & CPU time & $E_1$ &  $E_2$  \\
				& & (one time step) &  &  \\
				\hline
				1 & 0.8 & 37.2347 & 4.8130e-01 & 4.0692e-01 \\
				\hline
				2  & 0.8 & 63.5058 & 9.6106e-01 & 6.2141e-01\\
				\hline
				3 & 0.8 & 63.2127 & 2.3673e-01 & 2.4150e-01\\
				\hline
			\end{tabular} 
		\end{center}
		\caption{CPU runtime and errors in the concentration profile in 3D, B-char, final time $T=8$.}\label{tab:3DErr-ball}
	\end{table} 
	
	Here, we note that although the simulations were only performed over a mesh with $16\times16\times16$ cubes, with a time step $\delta t = 0.8$, the errors obtained in the numerical simulations in 3D is of a similar magnitude as those obtained in 2D. As with the tests in 2D, by using a continuous initial condition in the third test case, the error dropped by at least a factor 2.5 compared to that of the second test case. Moreover, the computational times, which ranged from 37 to 64 seconds for the three test cases, indicate that through the B-char method, characteristic-based schemes are not too costly to implement in 3D. 

These tests demonstrate that the B-char approach is a cost-effective and accurate way of implementing characteristic-based schemes such as the ELLAM, in both 2D and 3D.

	\begin{remark}[Number of balls used for B-char] 
Upon increasing the number of balls used for the B-char method, the initial approximation of the volumes becomes better. For the tests presented above, this allows us to achieve results with a similar quality even if we take $N<10$ iterations of \eqref{eq:volAdj} before solving the optimisation problem \eqref{eq:opti_masscons}. The choice here is between two options: increasing the number of balls and ball intersections, or increasing the number of iterations of \eqref{eq:volAdj}; we found that, for the tests presented above, the choice of 4 balls (in 2D) or 8 balls (in 3D) was producing a good accuracy/cost ratio. 
	\end{remark}

	\section{Conclusion and possible outlooks}
	
	In this paper, we have developed the B-char method, which is a cheap and efficient way to implement characteristic-based schemes, whilst preserving the important properties of local and global mass conservation, applicable in both 2D and 3D. Numerical tests were provided, and showed that the quality of the solutions obtained from B-char ELLAM is similar to the quality obtained by tracking and computing intersections of polygonal cells. The computational cost of the B-char method, however, is considerably lower than the computational cost of implementing the polygonal ELLAM. The tests run on benchmarks from the literature also showed that, in most situations, the B-char ELLAM produces rather accurate solutions, comparable to those obtained via other numerical schemes and much larger time steps.

Extension of the B-char method onto generic meshes can also be achieved, provided that we can design an algorithm to pack the balls inside the cells. Although the numerical tests were only performed on ELLAM type schemes, it would be interesting to explore the extension of the B-char method onto other characteristic-based schemes which are globally mass conserving, such as the MMOCAA. Future work would involve extending the B-char method onto non divergence-free velocity fields, and also those that are approximated by numerical schemes, e.g. piecewise polynomial velocity fields. The problem here lies with the fact that the divergence of the velocity field may be different from cell to cell, which makes it difficult to obtain a good approximation for the measure of the trace-back balls $|\widehat{B}_{K,s}|$. This requires further development and study of the idea proposed for approximating $\widehat{r}_{K,s}$ in \eqref{eq:rTrackedBall} for non-constant $\phi$, so that \eqref{eq:ball_int} gives a good enough initial approximation for the volume of the intersecting region $|F_{-\dtDisc}(K)\cap M|_\phi|$.

	\bibliographystyle{abbrv}
	\bibliography{ELLAM_Ball_Approx}
\end{document}